\documentclass[12pt,a4paper,oneside,fleqn]{article}
\usepackage{latexsym}
\usepackage{amstext,amsthm,amssymb}
\usepackage{amsmath}
\usepackage{amsfonts}
\usepackage[german,american]{babel}

\newcommand{\fa}{\qquad \forall \quad}
\newcommand{\lto}{\longrightarrow}

\newcommand{\eps}{\varepsilon}
\newcommand{\ro}{\varrho}
\newcommand{\fhi}{\varphi}
\newcommand{\teta}{\vartheta}

\newcommand{\tdelta}{\tilde{\delta}}
\newcommand{\tro}{\tilde{\ro}}
\newcommand{\tsigma}{\tilde{\sigma}}
\newcommand{\Wtilde}{\widetilde{W}}

\newcommand{\id}{{\rm id}}
\newcommand{\im}{{\rm Im}\,}

\newcommand{\stv}{\,\cdot\,}
\newcommand{\tp}{{\rm T}}

\newcommand{\heq}{\hspace{-2mm}&=&\hspace{-2mm}}

\newcommand{\hleq}{\hspace{-2mm}&\leq&\hspace{-2mm}}
\newcommand{\hlto}{\hspace{-2mm}&\lto&\hspace{-2mm}}

\newcommand{\hmapsto}{\hspace{-2mm}&\mapsto&\hspace{-2mm}}

\newcommand{\folgt}{\;\Rightarrow\;}
\newcommand{\gdw}{\;\Leftrightarrow\;}

\newcommand{\skal}[2]{\left\langle\, #1 \, ,\, #2 \,\right\rangle}
\newcommand{\poisson}[2]{\left\{ #1 \, ,\, #2 \right\}}

\newcommand{\maatrix}[4]{\left( \begin{array}{cc} #1 & #2 \\ #3 & #4
  \end{array} \right)}
\newcommand{\komp}[3]{#2 \stackrel{#1,\,{\rm\scriptscriptstyle
      compact}}{=\!=\!=\!\Longrightarrow} #3}

\newcommand{\kooomp}[3]{#2 \stackrel{#1,\,{\rm\scriptscriptstyle
      compact}}{=\!=\!=\!=\!=\!\Longrightarrow} #3}
\newcommand{\dis}{\displaystyle}

\newcommand {\B}{\mathcal{B}}
\newcommand {\C}{\mathbb C}

\newcommand {\D}{\mathcal{D}}

\newcommand {\M}{\mathcal{M}}
\newcommand {\N}{\mathbb N}
\newcommand {\Oo}{\mathcal{O}}
\newcommand {\Pp}{\mathcal{P}}
\newcommand {\Q}{\mathbb Q}
\newcommand {\R}{\mathbb R}
\newcommand {\Rr}{\mathcal{R}}
\newcommand {\Ss}{\mathcal{S}}

\newcommand {\U}{\mathcal{U}}
\newcommand {\V}{\mathcal{V}}

\newcommand {\Z}{\mathbb Z}

\newtheorem{theorem}{Theorem}[section]
\newtheorem{lemma}[theorem]{Lemma}
\newtheorem{corollary}[theorem]{Corollary}
\newtheorem{definition}[theorem]{Definition}
\newtheorem{remark}[theorem]{Remark}

\begin{document}
\title{An analytic KAM-Theorem}
\author{Joachim Albrecht}
\maketitle
\begin{abstract}
We prove an analytic KAM-theorem, which is used in \cite{albrecht07}, where the differential part of KAM-theory is discussed. Related theorems on analytic KAM-theory exist in the literature (e. g., among many others, \cite{moser69}, \cite{poeschel80}, \cite{salamon04}). The aim of the theorem presented here is to provide exactly the estimates needed in \cite{albrecht07}.
\end{abstract}
\tableofcontents\newpage
\renewcommand{\theequation}{\arabic{section}.\arabic{equation}}
\section{Formulation of the main theorem} \label{subak}
We consider Hamiltonian systems of the form
\begin{equation} \label{e1}
\dot{x} = H_y, \quad \dot{y} = -H_x.
\end{equation}
Here $x = (x_1, \ldots, x_n)$, $y = (y_1, \ldots, y_n)$, $\dot{x}$ and
$\dot{y}$ are vectors in $\R^n$ ($n \geq 2$) and $H=H(x,y)$ is a function from
$\R^{2n}$ to $\R$. 
We try to prove the existence of solutions of a system (\ref{e1}) under the
assumption, that it can be written as a sum
$ H = N + \widetilde{R}$
with a function 
\[ N(x,y) = a + \skal{\omega}{y} + \frac{1}{2} \skal{y\, Q(x)}{y} + \Oo(|y|^3),\]
\[( a\in \R, \, \omega \in \R^n, Q(x) \in \R^{n \times n},\, \skal{x}{y} := x_1 y_1 + \ldots + x_n y_n)\]
which we call normal form, and a remainder $\widetilde{R}$.
The dynamics of $N$ read
\[ \dot{x} = N_y = \omega + \Oo(|y|), \quad \dot{y} = -N_x = \Oo(|y|),\]
and are solved by
\[ t \mapsto (\omega t + \mbox{const.},0).\]
In case the frequencies $\omega_1, \ldots, \omega_n$ are rationally independent, such a solution is called quasi-periodic and it covers the Torus $\R^n/(2\pi\Z^n) \times \{0\}$ densely. 
KAM-Theory provides the means to prove, that many quasiperiodic
solutions survive the perturbation of the Hamiltonian. In our notation, the perturbed Hamiltonian is given by
\begin{equation*}
H(x,y) = a + \skal{\omega}{y} + \frac{1}{2} \skal{y \,Q(x)}{y} + R(x,y),
\end{equation*}
where $R$ denotes the sum of the terms of higher order of $N$ and the remainder $\widetilde{R}$. We prove the existence of quasiperiodic solutions of (\ref{e1}) for Hamiltonians of this kind.
\subsubsection*{Notations and Definitions}
For vectors $z = \left( z_1, \ldots,  z_\ell \right) \in \C^\ell$ we use the
$\ell_\infty$-norm \label{normvektor}
$|z| := \max_{1 \leq i \leq \ell} |z_i|$.
For matrices $Q = \left(q_{ij}\right)\in \C^{k \times \ell}$
we use the row-sum norm
\begin{equation*} 
|Q| := \max_{1 \leq i \leq k} \sum_{j=1}^\ell |q_{ij}|.
\end{equation*}
For arbitrary matrices $Q \in \C^{k \times \ell}$ and $P \in \C^{\ell \times m}$ the inequality $\left| Q \, P \right| \leq |Q| \, |P|$
holds. Transposed vectors and matrices are denoted with a superscript
``T''. \label{transp} For transposed matrices we have the estimate
$\left| Q^\tp \right| \leq k |Q|$, in which $Q$ has $k$ rows.
The product of two vectors $x$, $y \in \C^\ell$ is defined by
\[
\skal{x}{y} := \sum_{j=1}^\ell x_j \,y_j.
\]
Then we have $| \skal{x}{y} | \leq \ell\, |x|\,|y|$.
For the product of a vector $x \in \C^\ell$ and
a matrix $Q \in \C^{k \times \ell}$ the estimate
$| x\, Q^\tp| \leq |x|\,|Q|$
holds. Finally we have
\begin{equation} \label{opnorm}
|Q| = \max_{|z| \leq 1} |z Q^{\rm T}|.
\end{equation}
{\bf Domains and functions.}
\begin{definition}{\em \label{defbereiche}
Let $r$ and $s$ be positive numbers. We define
\[ \D(r,s) := \left\{ z=(x,y)\in \C^{2n}\,\left|\, |\im x|<r, |y|<s \right.\right\},\]
\[ \Ss(r) := \left\{ x\in \C^n\,\left|\, |\im x|<r\right.\right\},\]
\[ \Ss'(r) := \left\{ z\in \C^{2n}\,\left|\, |\im z|<r\right.\right\}.\]
Let $\Pp_m(r,s)$ be the set of all functions
\[f:\D(r,s)\lto \C^m, \quad z=(x,y) \mapsto f(z),\]
which are analytic, map real vectors to real values, and have period
$2\pi$ in the variables $x_1, \ldots, x_n$.\\
The set of all functions $f:\Ss(r)\to \C^m$, which are analytic, map real
vectors to real values and have period $2\pi$ in every variable, is
denoted by $\Pp_m(r)$.\\
The set of all functions $f:\Ss'(r)\to \C^m$, which are analytic, map real
vectors to real values and have period $2\pi$ in every variable, is
denoted by $\Pp_m'(r)$.\\\indent
 The definition shall hold for $m = n\times n$ as well. In case $m=1$
 we write $\Pp(r,s) := \Pp_1(r,s)$, $\Pp(r) := \Pp_1(r)$, and $\Pp'(r) := \Pp_1'(r).$
}\end{definition}
We denote the restriction of a function $f$ to a subset $\M$ of its domain with $\left. f\right|_\M$.\label{einschraenkung}
Notation of derivatives. Derivatives are denoted with a subscript, for example
\[ f_{x_1} = \frac{\partial f}{\partial x_1}, \quad f_x = (f_{x_1}, f_{x_2}, \ldots, f_{x_n}).\]
Hence, for a function $f = (f_1, \ldots, f_m) \in \Pp_m(r)$, $f_x$ is the  Jacobian. \label{ableitungen}
Finally we write for functions $t \mapsto (x_1(t), \ldots, x_n(t))$ depending
 on a single variable only
\[ \frac{dx}{dt} = \dot{x} = (\dot{x}_1, \ldots, \dot{x}_n) = (x_{1t}, \ldots, x_{nt}).\]
By our definition of the Jacobian we have $\dot{x} = x_t^{\rm
  T}$.\label{zeitabl}\\\\
{\bf Frequency vectors.}
The vector $\omega = (\omega_1, \ldots, \omega_n) \in \R^n$, which comes into
play as the first derivative of the Hamiltonian, is called {\em frequency
  vector}. To prove theorem \ref{hauptan} one has to assume that it satisfies
a sequence of {\em Diophantine inequalities}. That means, it has to be an
element of a set of the following type:
\begin{definition}{\em \label{omgamtau}
For $n\geq 2$, $\tau >0$, and $\gamma >0$ let
\[ \Omega(\gamma, \tau) := \left\{ \omega \in \R^n\,\bigg|\, |\skal{\omega}{k}| \geq \frac{\gamma}{|k|^\tau}\fa k \in \Z^n\setminus\{0\}\right\}.\]
}\end{definition}
\begin{remark}{\em \label{bemgamtau}
The following assertions hold (see \cite{ruessmann75} and the literature
given there):
\begin{enumerate}
\item In case $0 < \tau < n-1$, all sets $\Omega(\gamma, \tau)$, $\gamma>0$,
  are empty.
\item In case $\tau = n-1$, the $n$-dimensional Lebesgue measure of the set
  $\Omega(n-1) := \cup_{\gamma>0} \Omega(\gamma, n-1)$ is $0$. 
However, the intersection of every open subset of $\R^n$ with $\Omega (n-1)$
  has the cardinality of $\R$.
\item In case $\tau > n-1$, there exists a $\gamma = \gamma(\omega)>0$ with
  $\omega \in \Omega (\gamma, \tau)$ for almost every $\omega \in \R^n$.
\end{enumerate}
}\end{remark}\noindent
{\bf Simple canonical transformations}
\begin{definition}{\em \label{defsymplek}
Let $\U$ and $\V \subseteq \C^n$ be open connected sets. Let $J$ be the matrix
\begin{equation*} 
J=\left( \!\!\begin{array}{rc} 0 & E_n \\ -E_n & 0 \end{array}\! \right) \in
\C^{2n\times 2n}, \quad E_n \mbox{ the } (n\times n) \mbox{ identity matrix.}
\end{equation*}
We call a differentiable map
\[ Z: \U\times \V \longrightarrow \C^{2n}, \quad \zeta = (\xi,\eta) \mapsto z = Z(\zeta) \]
{\em symplectic transformation}, if for all $\zeta$ in $\U
\times \V$ the equation
\begin{equation} \label{defsympl}
Z_\zeta(\zeta)^{\rm T} \cdot J \cdot Z_\zeta(\zeta) = J
\end{equation}
holds.
}\end{definition}
\begin{definition}{\em \label{einfachetrafo}
Let $\U$, $\V \subseteq \C^n$ be open connected sets. We call an analytic symplectic transformation
\[ Z:\U\times \V \longrightarrow \C^{2n}, \quad \zeta = (\xi, \eta) \mapsto z = (x,y) = Z(\zeta) = (X(\zeta), Y(\zeta))\]
{\em simple canonical transformation}, if the map $\zeta=(\xi,\eta)
\mapsto X(\zeta)$ does not depend on $\eta$, which means $X=X(\xi)$.
}\end{definition}
Whenever the composition of two simple canonical transformations $Z_1$ and
$Z_2$ is possible,  $Z_1
\circ Z_2$ is a simple canonical transformation as well. If $Z_1$ and $Z_2$
have the property, that $(\xi,\eta) \mapsto Z_i(\xi, \eta) - (\xi, 0)$ has the period $2\pi$ in
$\xi_1, \ldots, \xi_n$ ($i =1,2$), so has
$(\xi,\eta) \mapsto Z_2 \circ Z_1 (\xi,\eta) - (\xi, 0)$.
\begin{theorem} \label{hauptan} {\bf Analytic KAM-theorem.}
\addcontentsline{toc}{subsubsection}{\hspace{1.325cm}Theorem \ref{hauptan}
  (analytic KAM-theorem)}
Let $\tau \geq n-1 \geq 1$, $\gamma > 0$, and $0 < s \leq r^{\tau+1} \leq
  1$. We consider the Hamiltonian $H \in \Pp(r,s)$, 
\begin{equation} \label{hami}
H(x,y) = a + \langle\omega, y\rangle + \frac{1}{2} \langle y\cdot Q(x) , y\rangle + R(x,y),
\end{equation}
where $a\in \R$, $\omega \in \Omega(\gamma,\tau)$, $Q \in \Pp_{n\times n}(r)$,
and $R\in \Pp(r,s)$. Let $C \in \R^{n\times n}$ be a non-singular matrix with 
\begin{equation} \label{nichtdegan}
\left|Q - C\right|_{\Ss(r)} \leq \frac{1}{4|C^{-1}|}.
\end{equation}
Then there exist positive constants $c_1$, $c_2, \ldots,$ $c_5$ depending on
$n$, $\tau$, $\gamma$, and $C$ only, such that for all $\teta$, $0 < \teta \leq c_1$, and
\begin{equation} \label{Man} M:=| R |_{\D(r, s)} \leq c_2 s^2 \teta
\end{equation}
the following holds: There exists a simple canonical transformation 
\[ W = (U, V): \D(r/2, s/2) \longrightarrow \D(r, s), \quad W-\id \in \Pp_{2n}(r/2, s/2)\]
with the estimate
\begin{equation} \label{trafo}
| W_\zeta - E_{2n}|_{\D(r/2, s/2)} \leq c_3 \teta.
\end{equation}
The transformed Hamiltonian $H_+ := H\circ W$ is an element of $\Pp(r/2, s/2)$
and has the form 
\begin{equation} \label{neu} 
H_+(\xi, \eta) = a_+ + \langle\omega, \eta\rangle + \frac{1}{2} \langle \eta \cdot Q_+(\xi) , \eta\rangle + R^\ast(\xi, \eta),
\end{equation}
where $a_+ \in \R$, $Q_+ \in \Pp_{n\times n}(r/2)$, and $R^\ast \in \Pp(r/2,
s/2)$. The functions $Q_+$ and $R^\ast$ fulfill the estimates
\begin{equation} \label{hesse}
\left| Q_+  - Q\right|_{\Ss(r/2)} \leq c_4 \teta,
\end{equation}
\begin{equation} \label{tayl3}
|R^\ast(\xi, \eta)| \leq c_5 M \frac{|\eta|^3}{s^3} \quad \mbox{for all } (\xi, \eta) \in \D(r/2, s/2).
\end{equation}
\end{theorem}
Assertion (\ref{tayl3}) means, that we can find solutions to the canonical
equations given by the Hamiltonian $H_+ = H\circ W$,
\begin{equation} \label{aha1}
\dot{\xi} = H_{+\eta},\quad
\dot{\eta} = -H_{+\xi}.
\end{equation}
Indeed, using the Landau symbol $\mathcal{O}$ we have $R^\ast =
\mathcal{O}(|\eta|^3)$, therefore (\ref{neu}) is the Taylor expansion of
$H_+$. So the equations (\ref{aha1}) can be written like this:
\[
\dot{\xi} = \omega + \Oo(|\eta|),\quad
\dot{\eta} = \Oo(|\eta|^2).
\]
We find the solution $\eta = 0,$ $\xi = \omega t + {\rm const.}$ 
It can be used to find a solution for the canonical equations corresponding to
the original Hamiltonian $H$,
\[
\dot{x} = H_y,\quad
\dot{y} = -H_x.
\]
Namely, the solution is $W(\xi, \eta) = W(\omega t+{\rm const.},0)$.\\\indent
The trick of theorem
\ref{hauptan}
is to get $\teta$ 
independent of $s$ in the estimates
(\ref{trafo}) and (\ref{hesse}).
This is essential to apply the theorem in differential KAM-theory.\\\indent
The fact, that $\omega$ can be kept fixed, is due to assumption
(\ref{nichtdegan}), for it causes $Q$ to be non-singular.
\setcounter{equation}{0}
\section{Motivation of the linearized equation}\label{sskaman}
We prove theorem
\ref{hauptan}
with Newton's method,
for its rapid convergence overcomes the influence of the so-called small divisors, see remarks
  \ref{klnennerbem} (page \pageref{klnennerbem}), \ref{klnennermu} (page
  \pageref{klnennermu}), and \ref{klnennermu2} (page \pageref{klnennermu2}).
To this end we have to establish a suitable linearised equation, which we now motivate. We write the Hamiltonian
(\ref{hami})
as a sum  
\[ H = N+R.\]
The summands are the
{\em normal form} 
\[ N(x,y) = a + \skal{\omega}{y} + \Oo(|y|^2),\]
and the -- small -- remainder
$R(x,y)$.
We have to find a sequence
$\left( Z_k \right)_{k\in \N}$
of symplectic transformations, such that the remainder gets smaller after every transformation. Write for $k\in \N_0$
\[ H=H_0, \quad H_k = N_k + R_k, \quad H_{k+1} := H_k \circ Z_{k+1},\]
where
$N_k$
again is a normal form (with
$a_k$ instead of $a$ and with the same 
$\omega$), and
$R_k$
is the remainder after the $k$-th step. When we set
\[ W_k := Z_1 \circ \ldots \circ Z_k, \quad W_0 := \id \quad (k\in \N),\]
we get
$H_k = H\circ W_k= N_k + R_k$. In case the limits
\[ R_k \longrightarrow 0, \quad W_k \lto W_\infty, \quad N_k \lto N_\infty \quad (k \to \infty)\]
exist with some symplectic transformation
$W_\infty$
and normal form
$N_\infty$, 
\[ H\circ W_\infty = N_\infty\]
follows and we are successful. In other words, we look for a root of the function
\[ \Rr(W,N) := H\circ W - N,\]
which is given by a pair of functions
$(W, N)$.
According to the above considerations, we try to find this root as a limit
\[ (W_\infty, N_\infty) = \lim_{k\to \infty} (W_k, N_k).\]
This leads to the problem to improve an approximate solution
$(W_k, N_k)$
to a better approximate solution
$(W_{k+1}, N_{k+1})$. For
$k \in \N_0$ we set\\[-2mm]
\parbox{11cm}{
\begin{eqnarray*}
W := W_k, \hspace{2.88cm} && N := N_k, \\ W_+ = W +\Delta W := W_{k+1} , && N_+ = N+ \Delta N:=N_{k+1},
\end{eqnarray*}}
\hfill
\parbox{3cm}{
\begin{equation} \label{dropk}
\end{equation}}\\
and obtain the new remainder as
\begin{eqnarray*}
\Rr(W_+, N_+) \heq H\circ (W+\Delta W) - N - \Delta N \\\heq \Rr(W,N) + H_z(W)\Delta W - \Delta N + {}\mbox{terms of higher order.}
\end{eqnarray*}
Linearisation means to solve the equation
\begin{equation} \label{lingl1}
\Rr(W,N) + H_z(W)\Delta W- \Delta N = 0.
\end{equation}
However, due to the term
$H_z(W) \Delta W$
 this is not possible in general. We have to separate further terms of higher order to get (\ref{lingl1}) simple enough. --
The following considerations are a simplified version of the approach presented in
\cite{ruessmann01}. (The situation in
\cite{ruessmann01}
is more complicated than the situation here because in \cite{ruessmann01}
the assumption
(\ref{nichtdegan})
is avoided.) We construct the symplectic transformations as flows of certain Hamiltonian systems. So we work with a function
$\Delta S = \Delta S(x,y)$
and consider the Hamiltonian system
\begin{equation} \label{konsys}
\dot{x} = \Delta S_y, \quad \dot{y} = - \Delta S_x.
\end{equation}
The solution of the respective initial value problem is denoted with
\[ z = (x,y) = (X(t, \xi, \eta), Y(t, \xi, \eta)) = Z(t, \xi, \eta), \quad Z(0, \xi, \eta) = (\xi, \eta) = \zeta.\]
Then, $t$ fixed, provided existence, the map
$\zeta \mapsto Z(t, \zeta)$
is a symplectic transformation (see appendix
\ref{gct}).
\begin{definition}{\em \label{pklammer}
Let
$f, g \in \Pp(r,s)$ or $f, g:\R^{2n}\to \R$
be differentiable functions. Then we define the {\em Poisson bracket} of $f$ and $g$ by
\[ \poisson{f}{g} := \skal{f_x}{g_y} - \skal{f_y}{g_x}. \]
}\end{definition}
For the moment let $F$ be a real valued, differentiable function. Then using (\ref{konsys})
we can replace a derivative with respect to time by a Poisson bracket as follows:
\begin{eqnarray} \nonumber
\frac{d}{d t} F(Z(t, \zeta)) \heq \skal{F_z (Z(t, \zeta))}{Z_t (t, \zeta)} \\ \nonumber
\heq \skal{F_x (Z(t,\zeta))}{X_t(t,\zeta)} + \skal{F_y(Z(t,\zeta))}{Y_t(t,\zeta)} \\ \nonumber
\heq \skal{F_x (Z(t,\zeta))}{\Delta S_y(Z(t,\zeta))} - \skal{F_y(Z(t,\zeta))}{\Delta S_x (Z(t,\zeta))} \\ \label{poissonleitetab}
\heq \poisson{F}{\Delta S}(Z(t,\zeta)).
\end{eqnarray}
Now assume the existence of a map 
$\zeta \mapsto Z(t, \zeta)$ for all $ 0 \leq t \leq 1$
and a set of allowed $\zeta$. The new transformation
$W_+ = W+\Delta W$
(see (\ref{dropk}))
shall be given by
$W_+(\zeta) := W(Z(1,\zeta))$.
$W$ being a symplectic transformation, 
$W_+$ will be a symplectic transformation as well.
With (\ref{poissonleitetab}) we get for 
$\Delta W$
the equation
\begin{eqnarray} \nonumber
\Delta W(\zeta) \heq W_+(\zeta) - W(\zeta) = W(Z(1,\zeta)) - W(\zeta) =
\int_0^1 \frac{d}{d t} W(Z(t,\zeta))\, dt \\ \label{deltawintegral}
\heq \int_0^1 \left( \poisson{W_{1}}{\Delta S}, \ldots, \poisson{W_{2n}}{\Delta S}\right) (Z(t,\zeta))\,dt.
\end{eqnarray}
Let us calculate
$\Rr(W_+, N_+)$ once more using (\ref{poissonleitetab}).
\begin{eqnarray*}
\Rr(W_+, N_+)(\zeta) \heq H \circ W_+ (\zeta) - N_+ (\zeta) = H\circ W(Z(1,\zeta)) - N(\zeta) - \Delta N(\zeta) \\
\heq \Rr(W,N) (Z(1,\zeta)) + N(Z(1,\zeta)) - N(\zeta) - \Delta N(\zeta) \\
\heq \big( \Rr(W,N) + \poisson{N}{\Delta S} - \Delta N \big) (\zeta) \\ &&{}+ \Rr(W,N)(Z(1,\zeta)) - \Rr(W,N)(\zeta)\\ &&{} + N(Z(1,\zeta)) - N(\zeta) - \left. \frac{d}{d t} N(Z(t,\zeta)) \right|_{t=0}.
\end{eqnarray*}
(The symbol
$\left.{}\right|_{t=0}$
means that the function has to be evaluated in the point
$t=0$\label{auswertung}.)
Like in
(\ref{deltawintegral})
we get
\begin{equation} \label{vereinfacher}
\Rr(W,N)(Z(1,\zeta)) - \Rr(W,N)(\zeta) = \int_0^1 \poisson{\Rr(W,N)}{\Delta
  S} (Z(t,\zeta))\,dt.
\end{equation}
Taylor's formula yields
\begin{equation} \label{vereinfachen}
 N(Z(1,\zeta)) - N(\zeta) - \left. \frac{d}{d t}
N(Z(t,\zeta)) \right|_{t=0} = \int_0^1 (1-t)
\frac{d^2}{d t^2} N(Z(t,\zeta))\,dt.
\end{equation}
When we use this, we obtain
\begin{eqnarray*}
\Rr(W_+, N_+)(\zeta) \heq \big( \Rr(W,N) + \poisson{N}{\Delta S} - \Delta N \big) (\zeta) + {}\\
&&{}+\int_0^1 \!\left( \poisson{\Rr(W,N)}{\Delta S} (Z(t, \zeta)) + (1-t)
\frac{d^2}{d t^2} N(Z(t,\zeta)) \right) dt.
\end{eqnarray*}
The time derivatives can be handled with
(\ref{poissonleitetab}),
\[ \frac{d^2}{d t^2} N(Z(t,\zeta)) = \frac{d}{d t} \poisson{N}{\Delta S} (Z(t,\zeta)) = \poisson{\poisson{N}{\Delta S}}{\Delta S} (Z(t,\zeta)),\]
\setlength{\mathindent}{0cm}\vspace{-.5cm}
\begin{eqnarray*} 
\Rightarrow \hspace{5mm} \Rr(W_+, N_+)(\zeta) \heq \big( \Rr(W,N) +
\poisson{N}{\Delta S} - \Delta N \big)(\zeta) +{} \\&&{}+ \int_0^1
\poisson{\Rr(W,N) + (1-t)\poisson{N}{\Delta S}}{\Delta S} (Z(t,\zeta))\,dt. 
\end{eqnarray*}
Hence we obtain the simplified linearised equation:
\setlength{\mathindent}{1cm}
\begin{equation} \fbox{$ \dis \Rr(W,N) + \poisson{N}{\Delta S} - \Delta N = 0$} \label{formlingl}
\end{equation}
This equation determines
$\Delta N$ and $\Delta S$.
Then $Z$ has to be calculated as the flow of
(\ref{konsys}).
This in turn determines
$W_+ = W\circ Z(1,\stv)$.
(\ref{formlingl}) being solved, the new remainder reads
\[ \Rr(W_+, N_+)(\zeta) = \int_0^1 \poisson{\Rr(W,N) + (1-t)\poisson{N}{\Delta S}}{\Delta S} (Z(t,\zeta))\,dt.\]
The inner Poisson bracket can be transformed with
(\ref{formlingl}),
for now 
\[ (1-t)\poisson{N}{\Delta S} = (1-t) \Delta N - (1-t) \Rr(W,N)\]
holds. So we can write
\begin{equation} \label{neuerrest}
\fbox{$\dis \Rr(W_+, N_+) (\zeta) = \int_0^1 \poisson{ t\,\Rr(W,N)+(1-t)\Delta N}{\Delta S}(Z(t,\zeta))\,dt$}
\end{equation} 
\setcounter{equation}{0}
\section{Solution of the linearized equation}
The solution of
(\ref{formlingl})
is based on the following theorem
\ref{linparttorus} 
from
\cite{ruessmann79c}
(in \cite{ruessmann79c} it is theorem 9.7).
\begin{definition}{\em \label{mittelwert}
Let
$r>0$ and $f: \Ss(r)\subseteq \C^n \to \C^m$, $x \mapsto f(x)$,
be a continuous function with period
$2\pi$ in $x_1, \ldots, x_n$.
We define the {\em mean}
$[f]$
of $f$ to be
\[ [f] := \left( \frac{1}{2\pi} \right)^n \int_0^{2\pi} \!\!\!\dots \int_0^{2\pi} f(x)\,dx_1\ldots dx_n.\]
}\end{definition}
\begin{theorem} \label{linparttorus}
Let
$\tau \geq n-1 \geq 1$, $\gamma>0$, $r>0$, $M>0$
and
$g:\Ss(r) \subseteq \C^n \to \C$
a
$2\pi$-periodic, analytic function with
$\left|g\right|_{\Ss(r)} \leq M$ and $[g]=0$. 
Let
$\omega \in \Omega(\gamma, \tau)$
(compare definition
\ref{omgamtau}).
Then there exists one and only one 
$2\pi$-periodic analytic function
$u:\Ss(r) \to \C$ 
with 
$[u]=0$ and 
\begin{equation} \label{klnennergl}
\skal{u_\xi}{\omega} = g.
\end{equation}
In addition there is a constant
$c_6 = c_6(n,\tau) >0$
with
\begin{equation} \label{klnennerabs}
\left|u\right|_{\Ss(r-\delta)} \leq \frac{c_6 M}{\gamma \delta^\tau}\fa
\delta \in (0, r).
\end{equation}
In case $g$ maps real vectors to real values, so does $u$.
\end{theorem}\noindent
\begin{remark} \label{klnennerbem}
Small divisors.
{\em Let us expand the given function $g$ and the solution $u$ into their Fourier series. These read, with coefficients
$ g_k$ 
and
$u_k \in \C$ $(k \in \Z^n \setminus \{0\})$, respectively,
\[ g(\xi) = \sum_{k\in \Z^n\setminus \{0\}}g_k e^{i\skal{k}{\xi}} \quad\mbox{and}
\quad u(\xi) = \sum_{k\in \Z^n\setminus \{0\}}u_k e^{i\skal{k}{\xi}}\fa \xi \in
\Ss(r).\]
The vanishing means of $g$ and $u$ amount to 
$g_0 = 0$ 
and
$u_0 = 0$,
respectively. The function $u$ can be differentiated term by term, so in
$\Ss(r)$
we get
\[ \skal{u_\xi(\xi)}{\omega} = \skal{ \sum_{k\in \Z^n\setminus \{0\}}i\,k\,u_k
  e^{i\skal{k}{\xi}}}{\omega} = \sum_{k\in \Z^n\setminus
  \{0\}}i\skal{k}{\omega}u_k e^{i\skal{k}{\xi}}.\]
Comparing coefficients with $g$ shows
$i\skal{k}{\omega}u_k=g_k$
for all
$k \in \Z^n\setminus\{0\}$. Hence
\begin{equation}  \label{klnennerreihe}
u(\xi) = \sum_{k\in \Z^n\setminus\{0\}} \frac{g_k}{i\skal{k}{\omega}}
  e^{i\skal{k}{\xi}}\fa \xi \in \Ss(r).
\end{equation}
So, if we took 
(\ref{klnennerreihe})
as an ansatz for the solution of the equation
$\skal{u_\xi}{\omega} = g$,
we had to proof convergence of this series. However, there is a serious obstacle: The divisors
 $i \skal{k}{\omega}$
 become very small -- in case the entries of 
 $\omega$ are not linear independent over $\Q$,
 there even exists some
$k \in \Z^n\setminus\{0\}$,
such that
$\skal{k}{\omega}$
vanishes: Therefore in this case there doesn't exist a 
$2\pi$-periodic analytic solution of
(\ref{klnennergl}).
\\
The meaning of theorem
\ref{linparttorus}
now is, that the series
(\ref{klnennerreihe})
indeed converges. The influence of the small divisors is represented by the factor
$c_6/(\gamma \delta^\tau)$
in estimate
(\ref{klnennerabs}).
} \end{remark}
\begin{theorem} \label{satzlingl}
Let
$\tau \geq n-1 \geq 1$, $\gamma >0$, $r>0$, $0<\delta < r/4$ and $0 <
s\leq \delta^{\tau+1} \leq 1$
be given. Suppose there is a constant
$M>0$
such that the function
$f\in \Pp(r,s)$
fulfills
\begin{equation} \label{N3vor}
\left| f \right|_{\D(r,s)} \leq M.
\end{equation}
Let
$N \in \Pp(r,s)$
be a function with
\begin{equation}
\label{N2vor} N(x,0) = N(0) \mbox{ and } N_y(x,0) = \omega \in
\Omega(\gamma,\tau)\fa x\in \Ss(r). 
\end{equation}
Finally, let 
$C \in \R^{n\times n}$
be a non-singular matrix with
\begin{equation} \label{N4vor} 
\left| N_{yy} - C \right|_{\D(r,s)} \leq \frac{1}{2 |C^{-1}|}.
\end{equation}
Then the equation
\begin{equation} \label{N*}
\fbox{$ \dis f+ \poisson{N}{\Delta S} - \Delta N=0$}
\end{equation}
possesses a solution, that is a pair of functions
$(\Delta S, \Delta N)$, with the properties:
\\
It is
$\Delta S (x,y)= \skal{\lambda}{x} +U(x) +\skal{V(x)}{y}$
with
$\lambda\in \R^n$ and $U\in \Pp(r)$, $V\in \Pp_n(r)$.
Especially the function
$(x,y) \mapsto \Delta S(x,y) - \skal{\lambda}{x}$ lies in $\Pp(r,s)$.
We have
$\Delta N \in \Pp(r,s)$,
\begin{equation}
\label{N2beh} \Delta N(x,0) = \Delta N(0) \mbox{ and } \Delta N_y(x,0) = 0
\fa x\in \Ss(r).
\end{equation}
There are constants
$c_7$, $c_8$, $\tilde{c}_9$, $c_{10}$ and $c_{11} > 0$,
such that the following estimates hold:
\begin{equation} \label{abssx}
\left|\Delta S_x\right|_{\D(r-4\delta, s)} \leq c_7\frac{M}{s},
\end{equation}
\begin{equation} \label{abssy}
\left|\Delta S_y\right|_{\Ss(r-3\delta)} \leq c_8\frac{M}{s \delta^\tau},
\end{equation}
\begin{equation} \label{absn0}
\left| \Delta N(0) \right| \leq \tilde{c}_9 \frac{M}{s},
\end{equation}
\begin{equation} \label{absn}
\dis \left| \Delta N - \Delta N(0)\right|_{\D(r-4\delta, s/2)} \leq c_{10} M,
\end{equation}
\begin{equation} \label{absnyy}
\dis \left| \Delta N_{yy} \right|_{\D(r-4\delta, s/4)} \leq c_{11} \frac{M}{s^2}.
\end{equation}
The constants
$c_j$ $(j \not= 9)$
only depend on
$n$, $\tau$, $\gamma$, and
$C$.
The constant
$\tilde{c}_9$
depends in addition on
$|\omega|$.
\end{theorem}\noindent
{\bf Proof.}
For
$\Delta S$
we make the ansatz
\begin{equation} \label{N++}
\Delta S(x,y) = \skal{\lambda}{x} + U(x) + \skal{V(x)}{y}.
\end{equation}
Here we try to obtain
$U \in \Pp(r)$ and $V\in \Pp_n(r)$
with
$[U] = 0$ and $[V] = 0$.
The vector
$\lambda \in \R^n$
has to be chosen suitable. We proceed in five steps.
\begin{enumerate}
\item Establish an equation to determine $U$.
\item Solve this equation.
\item Establish an equation to determine $V$.
\item Define $\lambda$ and solve the equation for $V$.
\item Define
$\Delta N$ 
and prove the properties of
$\Delta S$ and $\Delta N$.
\end{enumerate}
{\bf (1)} 
We deduce an equation for $U$. To this end we put $y=0$ in
(\ref{N*}). Assuming
$\Delta N(x,0) = \Delta N(0)$ for $x \in \Ss(r)$
(see (\ref{N2beh})) we obtain with
(\ref{N2vor})
\begin{eqnarray*}  
\lefteqn{f(x,0) + \poisson{N}{\Delta S}(x,0) - \Delta N(x,0) =} &&\\
\heq f(x,0) +
\skal{N_x}{\Delta S_y}(x,0) - \skal{N_y}{\Delta S_x}(x,0) - \Delta
N(0) \\ \heq f(x,0) - \skal{\Delta S_x(x,0)}{\omega} - \Delta N(0).
\end{eqnarray*}
This has to be zero. By
(\ref{N++})
that means for
$\Delta S$
\begin{equation} \label{Nvor}
f(x,0) - \skal{\lambda}{\omega} - \skal{U_x(x)}{\omega} - \Delta N(0) = 0.
\end{equation}
Well, with the help of theorem
\ref{linparttorus}
we can solve the equation
\begin{equation} \label{N+}
\skal{U_x(x)}{\omega} = f(x,0) - [f(\stv,0)].
\end{equation}
We take this equation to determine
$U$.\\
{\em 
Remark on the connection between equations
(\ref{Nvor}) and (\ref{N+}):}
Clearly
(\ref{Nvor}) and (\ref{N+})
are equivalent, if
\begin{equation} \label{ubedingt}
\Delta N(0) = \left[ f(\stv,0)\right] - \skal{\lambda}{\omega}.
\end{equation}
In step {\bf (4)} we will have to fix
$\lambda$ in such a way that the equation for $V$ is solvable, 
and then in step 
{\bf (5)}
define
$\Delta N$
such that
(\ref{ubedingt})
holds.
\\\\
{\bf (2)}
Solution of equation
(\ref{N+}).
The right hand side of 
(\ref{N+})
is bounded by
$2M$
because of (\ref{N3vor}).
Hence
Theorem
\ref{linparttorus}
yields a solution
$U\in \Pp(r)$
with
$[U]=0$ and 
\[ \left| U\right|_{\Ss(r-\delta)} \leq \frac{c_6 2M}{\gamma \delta^\tau}\fa
\delta \in (0,r).\]
With Cauchy's estimate (see lemma
\ref{cauchysabs} in the appendix)
we obtain
\begin{equation} \label{Npfund}
\left| U_x \right|_{\Ss(r-2\delta)} \leq \frac{2c_6 M}{\gamma
  \delta^{\tau+1}}\fa \delta \in (0, r/2).
\end{equation}
{\bf (3)} 
Now we have to find an equation for $V$. To this end we differentiate
(\ref{N*}) with respect to $y$ and put
$y=0$
to get
\begin{eqnarray*}
0 \heq f_y(x,0) + \poisson{N}{\Delta
    S}_y(x,0) - \Delta N_y(x,0)\\
\heq f_y(x,0) + 
    \skal{N_x}{\Delta S_y}_y(x,0) -
    \skal{N_y}{\Delta S_x}_y(x,0)
     - \Delta N_y(x,0) \\
\heq f_y(x,0) + \Delta S_y(x,0)\cdot N_{xy}(x,0) + N_x(x,0) \cdot \Delta S_{yy}(x,0) \\ && {} - \Delta S_x(x,0) \cdot N_{yy}(x,0) - N_y(x,0) \cdot \Delta S_{xy}(x,0) - \Delta N_y(x,0).
\end{eqnarray*}
The second summand vanishes because of
(\ref{N2vor}). The third summand is zero as well
by construction
(\ref{N++}). Therefore 
(\ref{N*})
implies
\begin{equation} \label{Numf}
f_y(x,0) - \Delta S_x(x,0) \cdot N_{yy}(x,0) - N_y(x,0) \cdot \Delta S_{xy}(x,0) - \Delta N_y(x,0)=0.
\end{equation}
Supposing
$\Delta N_y(x,0) = 0$ for $x\in \Ss(r)$
(compare (\ref{N2beh}))
we get with
(\ref{N2vor}) and (\ref{N++})
\[
f_y(x,0) - \left( \lambda + U_x(x)\right) \cdot N_{yy}(x,0) - \omega \cdot V_x^{\rm T}(x) = 0
\]
\vspace{-4mm}
\setlength{\mathindent}{0cm}
\begin{equation} \label{Nund}
\Leftrightarrow \hspace{4.5mm} \omega \cdot V_x^{\rm T}(x) = f_y(x,0) - \left( \lambda + U_x(x)\right) \cdot N_{yy}(x,0). 
\end{equation}
This is a system of $n$ equations which can be solved separately by theorem
\ref{linparttorus}, provided
\setlength{\mathindent}{1cm}
\begin{eqnarray*} 
0 \hspace{-2mm}&=&\hspace{-2mm} \left[ f_y(\stv,0) - ( \lambda + U_x)\cdot
  N_{yy}(\stv,0)\right]\\ 
\hspace{-2mm}&=&\hspace{-2mm} [f_y(\stv,0)] - [U_x\cdot N_{yy}(\stv,0)] - \lambda [N_{yy}(\stv,0)]
\end{eqnarray*} 
\vspace{-3mm}
\setlength{\mathindent}{0cm}
\begin{equation} \label{defvonlambda} 
\Leftrightarrow \hspace{4.5mm} 
\lambda [N_{yy}(\stv,0)] = [f_y(\stv,0)] - [U_x\cdot N_{yy}(\stv,0)].
\end{equation}
This equation has to be solved for $\lambda$.
\\\\
{\bf (4)} 
Definition of $\lambda$ and solution of
(\ref{Nund}). 
When $[N_{yy}(\stv,0)]$ is non-singular,
equation
(\ref{defvonlambda})
can be solved for $\lambda$.
We apply Lemma
\ref{invmatrix}
to
$[N_{yy}(\stv,0)]$.
By 
(\ref{N4vor})
\setlength{\mathindent}{1cm}
\[ \left| [N_{yy} (\stv,0)]-C \right| \leq \frac{1}{2 |C^{-1}|}\]
holds. So we can set  $S=C$, $P= [N_{yy}(\stv,0)]$, and
$h=1/2$ in the assumptions of
lemma
\ref{invmatrix}. It follows, that
$[N_{yy}(\stv,0)]^{-1}$ exists and that we have the estimate
\begin{equation} \label{NNyy}
\left| [N_{yy} (\stv,0)]^{-1}\right| \leq 2 |C^{-1}|.
\end{equation}
Therefore $\lambda$ can be defined as
\[ \lambda := \left( [f_y(\stv,0)] - [U_x \cdot N_{yy}(\stv,0)]\right) \cdot [N_{yy} (\stv,0)]^{-1}.\]
This choice guarantees, that the mean of the right hand side of
(\ref{Nund}) vanishes.
In order to apply theorem
\ref{linparttorus}
to
(\ref{Nund}),
we have to find an estimate for the right hand side of
(\ref{Nund}).
To begin with, 
(\ref{N3vor}) and Cauchy's estimate yield
\[ \left| f_y(\stv,0)\right|_{\Ss(r)} \leq \frac{M}{s}.\]
With respect to
$N_{yy}$ we observe
\[ 1 = \left| C C^{-1} \right| \leq |C||C^{-1}| \folgt \frac{1}{|C^{-1}|} \leq |C|,\]
hence with (\ref{N4vor}) we see
\begin{equation} \label{Ndollar}
\left| N_{yy} \right|_{\D(r,s)} \leq \left| N_{yy} -C\right|_{\D(r,s)} + |C| \leq \frac{1}{2 |C^{-1}|} + |C| \leq 2|C|.
\end{equation}
Together with
(\ref{Npfund}) and $s \leq \delta^{\tau+1}$
\begin{eqnarray} \nonumber
\left| f_y(\stv,0) - U_x\cdot N_{yy}(\stv,0)\right|_{\Ss(r-2\delta)} \hleq \frac{M}{s} + \frac{2 c_6 M}{\gamma \delta^{\tau+1}} 2 |C| \\\hleq \left( 1 + \frac{4 c_6 |C|}{\gamma} \right) \frac{M}{s} = c_{12} \frac{M}{s}\label{Nweigerlich} 
\end{eqnarray}
follows, where
\begin{equation} \label{defc12}
c_{12} := 1 + \frac{4 c_6 |C|}{\gamma}
\end{equation}
is a positive constant. This and
(\ref{NNyy}) give an estimate for $\lambda$, namely
\begin{equation} \label{Nweigerlambda}
|\lambda| \leq 2|C^{-1}|c_{12} \frac{M}{s}.
\end{equation}
The desired estimate for the right hand side of (\ref{Nund}) can be found using
(\ref{Ndollar}), (\ref{Nweigerlich}), and (\ref{Nweigerlambda}):
\begin{eqnarray} \label{Nproz} \nonumber
\lefteqn{\left| f_y(\stv,0) - (\lambda + U_x)\cdot N_{yy}(\stv,0)\right|_{\Ss(r-2\delta)} \leq \left| f_y(\stv,0) - U_x\cdot N_{yy}(\stv,0)\right|_{\Ss(r-2\delta)} +}\\& &{}+ |\lambda|\,\left| N_{yy}(\stv,0)\right|_{\Ss(r)} \leq c_{12} \frac{M}{s} + 4 |C|\,|C^{-1}| c_{12} \frac{M}{s}.
\end{eqnarray}
Now we can solve
(\ref{Nund}). Observe
\[ V(x) = (V_1(x),\ldots, V_n(x)), \quad \omega\cdot V_x^{\rm T}(x) = \left( \skal{\omega}{V_{1x}(x)}, \ldots, \skal{\omega}{V_{nx}(x)}\right).\]
Estimates for every
$V_i$ $(1\leq i \leq n)$
become estimates for $V$ for we use the maximum norm.
The right hand side of
(\ref{Nund})
is bounded on every substrip
$\Ss(r-\eps)$ of $\Ss(r)$
$(\eps \in (0, r))$,
because $f$, $U$, and $N$ are periodic in $x$. Therefore the solution $V$ exists on
$\Ss(r)$ and we have
$V\in \Pp_n(r)$
with the estimate
\begin{equation} \label{NabsV}
\left|V\right|_{\Ss(r-3\delta)} \leq \frac{c_6}{\gamma \delta^\tau} \left( c_{12} \frac{M}{s} + 4 |C|\,|C^{-1}| c_{12} \frac{M}{s} \right) = c_8 \frac{M}{s \delta^\tau}.
\end{equation}
Herein 
$c_8 = c_8(n, \tau, \gamma, C)$
is a positive constant. Further Cauchy's estimate yields
\begin{equation} \label{NabsVx}
\left| V_x \right|_{\Ss(r-4\delta)} \leq c_8 \frac{M}{s \delta^{\tau+1}}.
\end{equation}
{\bf (5)}
Now let us define
$\Delta S$
by
(\ref{N++}).
Then the assertions on the form of 
$\Delta S$
are fulfilled automatically. The definition
\[ \Delta N := f + \poisson{N}{\Delta S}\]
solves
(\ref{N*})
and 
$\Delta N \in \Pp(r,s)$
holds as well. Assertion
(\ref{N2beh})
is on the form of
$\Delta N$.
Using
(\ref{N2vor}), (\ref{N++}), and (\ref{N+})
we get
\begin{eqnarray} \nonumber 
\Delta N(x,0) \hspace{-2mm}&=&\hspace{-2mm} f(x,0) + \poisson{N}{\Delta
  S}(x,0) \\
\nonumber \heq f(x,0) + \skal{N_x}{\Delta S_y}(x,0) -
  \skal{N_y}{\Delta S_x}(x,0) \\
\nonumber \heq f(x,0) - \skal{\omega}{\Delta S_x(x,0)} \\
\nonumber \heq f(x,0) - \skal{\lambda}{\omega} - \skal{U_x(x)}{\omega} \\
\hspace{-2mm}&=&\hspace{-2mm} [f(\stv,0)] - \skal{\lambda}{\omega}. \label{heinzelmann}
\end{eqnarray}
This is obviously independent of $x$. So we may write
$\Delta N(x,0) = \Delta N(0)$
for all
$x\in \Ss(r)$.
Incidentally the calculation shows, that
(\ref{ubedingt}) is fulfilled and that
solving
(\ref{N+})
solves
(\ref{Nvor})
as well. -- In 
(\ref{Numf}) we have seen, that equation
(\ref{N*}), which we have proven in the meantime, implies
\[ \Delta N_y(x,0) = f_y(x,0) - \Delta S_x(x,0) \cdot N_{yy}(x,0) - N_y(x,0) \cdot \Delta S_{xy}(x,0).\]
Therefore
(\ref{N++}) and (\ref{Nund}) yield
\[ \Delta N_y(x,0) = f_y(x,0) - \left( \lambda + U_x(x) \right) \cdot N_{yy}(x,0) - \omega \cdot V_x^{\rm T}(x) = 0,\]
and
(\ref{N2beh}) is shown. We turn to the estimates for the derivatives of
$\Delta S$. By definition
(\ref{N++})
$\Delta S_y = V$,
so
(\ref{NabsV}) means
\[ \left| \Delta S_y \right|_{\Ss(r-3\delta)} \leq c_8 \frac{M}{s \delta^\tau}.\]
This is
(\ref{abssy}). We have
$\Delta S_x(x,y) = \lambda + U_x(x) + y\cdot V_x(x)$.
With
(\ref{Nweigerlambda}), (\ref{Npfund}), (\ref{NabsVx}) and the assumption $s \leq \delta^{\tau+1}$
we calculate
\begin{eqnarray*}
\left| \Delta S_x \right|_{\D(r-4\delta, s)} \hleq |\lambda| + \left| U_x \right|_{\Ss(r-2\delta)} + n s \left| V_x \right|_{\Ss(r-4\delta)} \\\hleq 2 |C^{-1}| c_{12} \frac{M}{s} + \frac{2 c_6 M}{\gamma s} + n s c_8 \frac{M}{s^2} = c_7 \frac{M}{s},
\end{eqnarray*}
where
$c_7 = c_7(n, \tau, \gamma, C)$
is a positive constant. This proves 
(\ref{abssx}). The estimates for
$\Delta N$ and $\Delta N_{yy}$
remain. In
(\ref{heinzelmann})
we have seen
$\Delta N(0) = [f(\stv,0)] - \skal{\lambda}{\omega}$.
According to
(\ref{N3vor}) and (\ref{Nweigerlambda})
this yields
\[ |\Delta N(0)| \leq M + 2 n|\omega| |C^{-1}| c_{12} \frac{M}{s} \leq \tilde{c}_9 \frac{M}{s},\]
where
$\tilde{c}_9 = \tilde{c}_9(n,\tau, \gamma, C, |\omega|)$
again is a positive constant. Hence 
(\ref{absn0}) holds. In order to show
(\ref{absn})
we use
(\ref{N++}), (\ref{N+}) and (\ref{heinzelmann})
to get
\begin{eqnarray*} 
\skal{\Delta S_x(x,y)}{\omega} \heq \skal{\lambda}{\omega} + \skal{U_x(x)}{\omega} + \skal{y\cdot V_x(x)}{\omega} \\ \heq \skal{\lambda}{\omega} + f(x,0) -[f(\stv,0)] +\skal{y}{\omega \cdot V_x^{\rm T}(x)} \\\heq f(x,0) + \skal{y}{\omega \cdot V_x^{\rm T}(x)} - \Delta N(0).
\end{eqnarray*}
With
(\ref{Nund}) and (\ref{Nproz})
we obtain
\begin{equation} \label{defc13strich}
\left| \skal{\Delta S_x}{\omega} + \Delta N(0)\right|_{\D(r-2\delta, s)} \leq
M + ns(c_{12} + 4 |C|\,|C^{-1}| c_{12})\frac{M}{s} = c_{13} M,
\end{equation}
where
\begin{equation} \label{defc13}
c_{13} := 1+ nc_{12}\left( 1+4|C|\,|C^{-1}|\right).
\end{equation}
Let us for the moment denote the function
$y \mapsto \skal{\omega}{y}$
by
$g_\omega$.
Then we can write
\begin{eqnarray} \nonumber
\Delta N \heq f+\poisson{N}{\Delta S} = f + \skal{N_x}{\Delta S_y} - \skal{N_y}{\Delta S_x}\\ \nonumber
\heq  f + \skal{ \left( N- g_\omega -N(0)\right)_x}{\Delta S_y} \\ \nonumber &&{}-
\skal{ ( N- g_\omega -N(0))_y}{\Delta S_x} - \skal{\omega}{\Delta S_x}\\ \label{dirkia}
\hspace{-2mm}&=&\hspace{-2mm} f + \poisson{N-g_\omega - N(0)}{\Delta S}-\skal{\omega}{\Delta S_x}.
\end{eqnarray}
Let us have a closer look at the first entry of the Poisson bracket. 
We can write
\begin{equation} \label{dirkib}
N(x,y) - \skal{\omega}{y} - N(0) = N(x,y) - \skal{N_y(x,0)}{y} - N(x,0) =: h(x,y)
\end{equation}
for all
$(x,y)\in \D(r,s)$
because of (\ref{N2vor}).
This defines a function
$h\in \Pp(r,s)$
with
$h(x,0) = 0$ and $h_y(x,0) = 0$
for all
$x\in \Ss(r)$.
Taylor's formula yields 
\[ |h(x,y)| \leq \left| \int_0^1 \frac{(1-\sigma)^2}{2} \skal{y \cdot
  h_{yy}(x,\sigma y)}{y}\,d\sigma \right| \leq \frac{1}{2} n|s|^2
\left|h_{yy}(x,\stv)\right|_{\{y\in \C^n\,|\, |y| <s\}}\]
for all  $(x,y) \in \D(r,s)$, from which we conclude with
(\ref{Ndollar})
\[ \left| h \right|_{\D(r,s)} \leq |C| n s^2.\]
Cauchy's estimate results in
\begin{equation} \label{dirkic}
\left| h_x \right|_{\D(r-\delta, s)} \leq \frac{|C| n s^2}{\delta}, \mbox{ and } \left| h_y \right|_{\D(r, s/2)} \leq 2 |C| n s.
\end{equation}
Now, 
(\ref{dirkia}) and (\ref{dirkib})
show
\begin{eqnarray*}
\Delta N - \Delta N(0) \heq f + \poisson{N-g_\omega - N(0)}{\Delta S} -
\skal{\omega}{\Delta S_x} - \Delta N(0) \\
\heq f + \poisson{h}{\Delta S} - \left( \skal{\omega}{\Delta S_x} + \Delta
N(0)\right) \\
\heq f + \skal{h_x}{\Delta S_y} - \skal{h_y}{\Delta S_x} - \left( \skal{\omega}{\Delta S_x} + \Delta N(0)\right).
\end{eqnarray*}
When we put the estimates for $f$, $\Delta S_y$ and $\Delta S_x$, 
(\ref{dirkic}), (\ref{defc13strich}), and (\ref{defc13}) together, we get
\begin{eqnarray*}
\left| \Delta N- \Delta N(0)\right|_{\D(r-4\delta, s/2)} \hspace{-2mm}&\leq&\hspace{-2mm} M + n\cdot\frac{|C| n s^2}{\delta} \cdot c_8 \frac{M}{s\delta^\tau} + n\cdot 2|C| ns \cdot c_7 \frac{M}{s} + c_{13} M\\ \hspace{-2mm}&\leq&\hspace{-2mm} c_{10} M,
\end{eqnarray*}
where
\begin{equation} \label{defc10}
c_{10} := 1 + n^2 |C| \left( 2c_7 +c_8\right) + c_{13}
\end{equation}
is a positive constant. This proves
(\ref{absn}).
Now
(\ref{absnyy})
is a consequence of Lemma
\ref{cauchysabs},
\[ \left| \Delta N_{yy} \right|_{\D(r-4\delta, s/4)} \leq \frac{8}{s}
\left| \Delta N_y \right|_{\D(r-4\delta, (3/8)s)} \leq 64 \frac{M c_{10}}{s^2}.\]
It remains only to set
$c_{11} = c_{11}(n,\tau,\gamma, C) := 64 c_{10} >0$
to finish the proof.
\hfill $\Box$
\setcounter{equation}{0}
\section{The inductive lemma}
In this section we construct a sequence of symplectic transformations and proceed in three steps. At first we prove theorem
\ref{induktionslemma}. It deals with a transformation $Z$, which transforms a given Hamiltonian $H$ into
$H_+ = H\circ Z$. Next we find sequences of numbers
$\left( r_k \right)$, $\left( \delta_k\right)$, $\left( s_k \right)$, and $\left( M_k\right)$,
such that theorem
\ref{induktionslemma}
can be applied repeatedly. That means that the obtained function
$H_+$ can be again inserted in the assumptions of theorem \ref{induktionslemma} as a new function $H$. The third step is to summarize the results and describe the inductive process for all $k\in \N_0$ in form of the inductive lemma 
\ref{induktionssatz}.
\begin{theorem} \label{induktionslemma}
Let
$\tau \geq n-1 \geq 1$, $\gamma>0$, $r >0$, $0<\delta < r/6$,
$0<s\leq \delta^{\tau+1} \leq 1$,
and
$0 < r_+ \leq r -6\delta$ and $0<s_+ \leq s/8$.
We consider a function
$H\in \Pp(r,s)$,
$H = N+R$
with
$N$, $R \in \Pp(r, s)$
and
\setlength{\mathindent}{1cm}
\begin{equation} \label{ind1vorA}
N(x,y) = a + \skal{\omega}{y} + \Oo(|y|^2), 
\end{equation}
where
$a\in \R$ and $\omega \in \Omega(\gamma, \tau)$ is assumed.
Further we assume the existence of a non-singular matrix
$C \in \R^{n\times n}$ with
\begin{equation} \label{ind1vorB}
\left| N_{yy} - C \right|_{\D(r, s)} \leq \frac{1}{2 |C^{-1}|}.
\end{equation}
The remainder $R$ has to be bounded by a constant
$M > 0$
with
\begin{equation} \label{ind1vorC}
\left| R \right|_{\D(r,s)} \leq M \leq \frac{1}{16} \frac{1}{ c_7 + c_8}\,s^2.
\end{equation}
Herein the constants
$c_7$ and $c_8$ are given by Theorem
\ref{satzlingl} (see (\ref{abssx}) and (\ref{abssy})).
Then there exists a simple canonical transformation 
(see definition \ref{einfachetrafo})
\begin{eqnarray} \label{trafoseigen}
Z:\D(r_+,s_+) \hlto \D(r-5\delta, s/4), \quad Z - \id \in \Pp_{2n}(r_+, s_+),\\\nonumber
\zeta = (\xi, \eta) \hmapsto Z(\xi, \eta), 
\end{eqnarray}
such that the transformed Hamiltonian
$H_+ = H\circ Z$
is an element of
$\Pp(r_+,s_+)$
and
$H_+ = N_+ + R_+$
holds, where
$N_+$, $R_+ \in \Pp(r_+, s_+)$, and
\begin{equation} \label{h+taylor}
N_+ (\xi, \eta) = a_+ + \skal{\omega}{\eta} + \Oo(|\eta|^2)
\end{equation}
with some
$a_+ \in \R$.
The following estimates hold:
\begin{equation} \label{ind1Z}
\left| Z_\zeta \right|_{\D(r_+, s_+)} \leq \exp \left( c_{14} \frac{M}{s^2} \right),
\end{equation}
\begin{equation} \label{ind1Z-E}
\left| Z_\zeta - E_{2n} \right|_{\D(r_+, s_+)} \leq c_{14} \frac{M}{s^2} \exp \left( c_{14} \frac{M}{s^2} \right),
\end{equation}
\begin{equation} \label{ind1a+-a}
\left| a_+ -a \right| \leq \widetilde{c}_9 \frac{M}{s},
\end{equation}
\begin{equation} \label{ind1Q+-Q}
\left| N_{+\eta\eta} - N_{\eta\eta} \right|_{\D(r_+, s_+)} \leq c_{11} \frac{M}{s^2},
\end{equation}
\begin{equation} \label{ind1R+}
\left| R_+ \right|_{\D(r_+,s_+)} \leq c_{15} \frac{M^2}{s^2}.
\end{equation}
The constants
$\widetilde{c}_9$ and $c_{11}$
are given by Theorem
\ref{satzlingl} (see (\ref{absn0}) and (\ref{absnyy})), and
$c_{14}$, $c_{15}$
are positive constants depending on
$n$, $\tau$, $\gamma$, and
$C$
only. Finally, if the partial derivatives $W_\xi$ and $W_\eta$ of the function
$W= W(\xi, \eta) : \D(r,s) \to \C^{2n}$
are continuous and bounded by 
$K_1>0$, then
$\Delta W:= W\circ Z - W$
satisfies
\begin{equation} \label{ind1DeltaW}
\left| \Delta W \right|_{\D(r_+, s_+)} \leq n K_1 (c_7+c_8) \frac{M}{s \delta^\tau}.
\end{equation}
\end{theorem}
\begin{remark}{\em \label{klnennermu}
We see the success of our approach in estimate
(\ref{ind1R+}), for the magnitude $M$ of the old remainder enters quadratically. This is due to Newton's method. The disturbing influence of the small divisors (compare remark \ref{klnennerbem})
is seen in the factor
$1/s^2$.
}\end{remark} \noindent
{\bf Proof of theorem \ref{induktionslemma}.}
We solve the linearized equation
\setlength{\mathindent}{1cm}
\begin{equation} \label{indlinear}
R + \poisson{N}{\Delta S} - \Delta N = 0
\end{equation}
by means of theorem
\ref{satzlingl}. Let us check the assumptions of that theorem. We apply the constants 
$\tau$, $\gamma$, $\delta$, $r$, $s$, and $M$ as they are in theorem
\ref{satzlingl}, such that the assumptions on those constants are fulfilled.
Further we insert
$f=R$ and $N= H-R$.
Now, 
$R$, $N \in \Pp(r,s)$
and from
(\ref{ind1vorA})
$N(x,0) = N(0)=a$ and $N_y(x,0) = \omega \in \Omega(\gamma,\tau)$ hold for all $x\in \Ss(r)$.
With
(\ref{ind1vorB}) and (\ref{ind1vorC})
all assumptions of theorem
\ref{satzlingl}
are met. Hence we obtain a solution
$(\Delta S, \Delta N)$
of
(\ref{indlinear})
with all the properties asserted in theorem
\ref{satzlingl},
especially the estimates
(\ref{abssx})
to
(\ref{absnyy}).\\
The construction of $Z$ proceeds like it is described in the appendix, see theorem
\ref{trafoterminator} in section \ref{gct}.
Theorem
\ref{trafoterminator} can be applied with
\begin{equation} \label{defc8}
K = (c_7 + c_8) \frac{M\delta}{s} >0,
\end{equation}
\[ \ro = r-4\delta,\,\sigma = s/4, \mbox{ and } F= \left. \Delta S
\right|_{\D(\ro, \sigma)} \in \Pp(\ro,\sigma).\]
We have $2\delta < \ro$ because of 
$\delta < r/6$ and
$0 < \sigma \leq \delta$
from
$0 < s \leq \delta^{\tau+1} \leq 1$.
(\ref{defc8}) and
(\ref{ind1vorC}) show
\[ \frac{\sigma \delta}{2 K} = \frac{\sigma \delta}{2}\cdot \frac{s}{(c_7 +
  c_8) M \delta} = \frac{s^2}{8 (c_7+c_8) M} \geq 2 > 1.\]
The function $F$ is affine linear in $y$, as is
$\Delta S$. 
We use
(\ref{abssx})
to get
\[ \left| F_x \right|_{\D(\ro,\sigma)} = \left| \Delta S_x \right|_{\D(\ro, \sigma)} \leq c_7 \frac{M}{s} \leq (c_7 + c_8) \frac{M \delta}{s} \cdot \frac{1}{\delta} = \frac{K}{\delta},\]
and
(\ref{abssy})
yields
\[ \left| F_y \right|_{\D(\ro,\sigma)} \leq \left| \Delta S_y \right|_{\Ss(r-3\delta)} \leq c_8 \frac{M}{s \delta^\tau} \leq (c_7 + c_8) \frac{M\delta}{\delta^{\tau+1}} \cdot \frac{1}{s} \leq \frac{K}{s} < \frac{4 K}{s} = \frac{K}{\sigma}.\]
So $F$ fulfills the assumptions
(\ref{vorh}) of theorem \ref{trafoterminator}, which can be applied now.
According to
(\ref{hamflussorig})
we obtain simple canonical transformations
\\
\begin{minipage}{10cm}
\[ Z(t,\stv):\D(r-6\delta, s/8) \lto \D(r-5\delta, s/4),\]\vspace{-7mm}
\[ Z(t,\stv) -\id \in
\Pp_{2n}(r-6\delta, s/8) \quad ( 0\leq t
< 2).\]
\end{minipage}
\begin{minipage}{5.5cm}
\begin{equation} \label{hamflusseinsatz}
\end{equation}
\end{minipage}\\[2mm]
With (\ref{defc8}) we calculate
\[ \frac{2nK}{\delta \sigma} = \frac{2\cdot 4n(c_7+c_8) M}{s^2} = c_{14} \frac{M}{s^2},\]
wherein $c_{14} = 8n(c_7+c_8)$ is a positive constant.
This can be put into the estimates
(\ref{flussab1})
and
(\ref{flussab2})
of theorem
\ref{trafoterminator}
to infer
\begin{equation} \label{ind1Zbew}
\left| Z_\zeta(t,\stv)\right|_{\D(r-6\delta, s/8)} \leq \exp\left( c_{14}
\frac{M}{s^2}\,t\right)\fa t \in [0, 2),
\end{equation}
\begin{equation} \label{ind1Z-Ebew}
\left| Z_\zeta(t,\stv)-E_{2n}\right|_{\D(r-6\delta, s/8)} \leq
c_{14}\frac{M}{s^2}\exp\left( c_{14} \frac{M}{s^2}\,t\right)\fa t
\in [0,1]
\end{equation}
for the maps given in
(\ref{hamflusseinsatz}).
Now we define 
$Z$
to be the function
$Z(1,\stv)$
restricted to
$\D(r_+, s_+)$.
Than $Z$ has the properties
(\ref{trafoseigen})
because of
(\ref{hamflusseinsatz}).
(\ref{ind1Zbew})
and
(\ref{ind1Z-Ebew})
cause $Z$ to meet the estimates
(\ref{ind1Z}) and (\ref{ind1Z-E}).\\
We set for all
$\zeta \in \D(r_+, s_+)$
\[H_+(\zeta) := (H\circ Z)(\zeta),\quad N_+(\zeta) := N(\zeta) + \Delta N(\zeta), \quad R_+(\zeta) := H_+(\zeta) - N_+(\zeta),\] 
(observe $N= H-R$). We deduce the properties of 
$N_+$ from the properties of
$\Delta N$
formulated in theorem
\ref{satzlingl}.
$\Delta N \in \Pp(r,s)$
implies
$N_+ \in \Pp(r_+, s_+)$. Furthermore,
\[ N_+(\xi,0) = N(\xi,0) + \Delta N(\xi,0) = a + \Delta N(0) =: a_+\fa \xi \in S(r_+).\]
(\ref{ind1a+-a}) is a consequence of
(\ref{absn0}):
\[ |a_+ - a| = \left| \Delta N(0)\right| \leq \widetilde{c}_9 \frac{M}{s}.\]
Next we see 
\[ N_{+y}(\xi,0) = N_y(\xi,0) + \Delta N_y(\xi,0) = \omega\fa \xi \in \Ss(r_+).\]
So the Taylor expansion of
$N_+$
is given by
\[ N_+( \xi,\eta) = a_+ + \skal{\omega}{\eta} + \Oo(|\eta|^2),\]
which is
(\ref{h+taylor}).
Estimate
(\ref{ind1Q+-Q})
follows from
(\ref{absnyy}):
\[ \left| N_{+\eta\eta} - N_{\eta\eta}\right|_{\D(r_+, s_+)} = \left| \Delta
N_{\eta\eta} \right|_{\D(r_+, s_+)} \leq c_{11}\frac{M}{s^2}.\]
Now we check
$R_+ \in \Pp(r_+,s_+)$: 
$R_+$ is an analytic function, which maps real vectors to real values, and
we have for all
$1 \leq j \leq n$
\begin{eqnarray*}
R_+(\xi + 2\pi e_j, \eta) \heq H(Z(\xi+2\pi e_j, \eta)) - N_+(\xi + 2\pi e_j, \eta) \\ \heq H(Z(\xi,\eta) + (2\pi e_j,0)) - N_+(\xi, \eta) \\ \heq H(Z(\xi, \eta)) - N_+ (\xi, \eta) = R_+(\xi,\eta),
\end{eqnarray*}
which is the desired periodicity. In order to prove
(\ref{ind1R+})
we recalculate 
(\ref{neuerrest}) -- we redo the calculations of section
\ref{sskaman}
with our functions, which are well-defined in the meantime, and use
(\ref{poissonleitetab}), (\ref{vereinfacher}), (\ref{vereinfachen}), and
(\ref{indlinear}):
\begin{eqnarray}
\nonumber R_+(\zeta) \heq H_+(\zeta) -N_+(\zeta) = H\circ Z(\zeta) -
N_+(\zeta) = H \circ Z(1,\zeta) - N(\zeta) - \Delta N(\zeta) \\
\nonumber \heq R(Z(1,\zeta)) + N(Z(1,\zeta)) - N(\zeta) - \Delta N(\zeta) \\
\nonumber \heq \left( R + \poisson{N}{\Delta S} - \Delta N\right)(\zeta) +
R(Z(1,\zeta)) - R (\zeta) \\
\nonumber &&{} + N(Z(1,\zeta)) - N(\zeta) - \left. \frac{d}{d t}
N(Z(t,\zeta))\right|_{t=0} \\
\nonumber \heq \int_0^1 \poisson{R}{\Delta S} (Z(t,\zeta))\,dt + \int_0^1
(1-t) \frac{d^2}{d t^2} N(Z(t,\zeta))\,dt\\
\nonumber \heq \int_0^1 \poisson{R + (1-t) \poisson{N}{\Delta S}}{\Delta
  S}(Z(t,\zeta))\,dt \\
\label{r+int} \heq \int_0^1 \poisson{tR+(1-t)\Delta N}{\Delta S}
(Z(t,\zeta))\,dt\fa \zeta \in \D(r_+, s_+).
\end{eqnarray}
To estimate the integrand we set for 
$t \in [0,1]$
\[ F_{(t)} := tR + (1-t) (\Delta N - \Delta N(0)) \in \Pp(r,s).\]
Then our assumption
(\ref{ind1vorC})
and
(\ref{absn}) 
lead to
\[ \left| F_{(t)} \right|_{\D(r-4\delta, s/2)} \leq tM + (1-t) c_{10} M \leq (1+c_{10}) M\fa t\in [0,1].\]
We use Cauchy's estimate to get for all
$t\in [0,1]$
\[ \left| F_{(t)x} \right|_{\D(r-5\delta, s/2)} \leq (1+c_{10}) \frac{M}{\delta}, \quad \left| F_{(t)y}\right|_{\D(r-4\delta, s/4)} \leq 4 (1+c_{10})\frac{M}{s}.\]
Together with
(\ref{abssx})
and
(\ref{abssy})
we obtain for all $t\in [0,1]$
\begin{eqnarray*}
\left| \poisson{F_{(t)}}{\Delta S} \right|_{\D(r-5\delta, s/4)} \hleq n \left( \left| F_{(t)x}\right|_{\D(r-5\delta, s/2)} \left| \Delta S_y\right|_{\Ss(r-3\delta)} +\right.\\ &&{}\left. +\left| F_{(t)y} \right|_{\D(r-4\delta, s/4)} \left| \Delta S_x\right|_{\D(r-4\delta, s)}\right) \\ \hleq n(1+c_{10}) \left( \frac{M}{\delta} c_8 \frac{M}{s \delta^\tau} + \frac{4M}{s}c_7 \frac{M}{s}\right) \leq c_{15} \frac{M^2}{s^2},
\end{eqnarray*}
where
\begin{equation} \label{defc15}
c_{15} := n(1+c_{10})(4c_7+c_8)
\end{equation}
is a positive constant.
Now,
\[ \poisson{tR + (1-t)\Delta N}{\Delta S} = \poisson{F_{(t)}}{\Delta S}\fa t \in [0,1],\]
and we have
$Z(t,\zeta) \in \D(r-5\delta, s/4)$
for all
$t\in [0,1]$ and $\zeta \in \D(r_+, s_+)$
by
(\ref{hamflusseinsatz}).
So we can deduce the estimate
(\ref{ind1R+}) for $R_+$ from
(\ref{r+int}).\\
Finally we have to show
(\ref{ind1DeltaW}).
The estimates for
$W_\xi$ and $W_\eta$ become estimates for
$W_{j\xi}$ and $W_{j\eta}$ $(1\leq j \leq 2n)$, because we use the row-sum norm.
Hence our assumptions read
\[ \left| W_{j\xi}\right|_{\D(r, s)} \leq K_1 \mbox{ and } \left|
W_{j\eta}\right|_{\D(r, s)} \leq K_1\fa 1 \leq j \leq 2n.\]
(\ref{deltawintegral})
implies 
for all $1\leq j \leq 2n$ and $\zeta \in \D(r_+, s_+)$ 
\[ \Delta W_j(\zeta) = \int_0^1 \poisson{W_j}{\Delta S}(Z(t,\zeta))\,dt.\]
So, writing
$\Delta S_\xi := \Delta S_x$ and $\Delta S_\eta := \Delta S_y$,
we obtain with
(\ref{abssx}) and (\ref{abssy})
\begin{eqnarray*} 
\left| \Delta W_j\right|_{\D(r_+, s_+)} \heq \left| \int_0^1
\poisson{W_j}{\Delta S} (Z(t,\stv))\,dt \right|_{\D(r_+, s_+)} \\
\hleq \int_0^1 \left| \poisson{W_j}{\Delta S}\right|_{\D(r-5\delta,
  s/4)}dt\\
\hleq \left| \skal{W_{j\xi}}{\Delta S_\eta}\right|_{\D(r-5\delta, s/4)} +
\left| \skal{W_{j\eta}}{\Delta S_\xi}\right|_{\D(r-5\delta, s/4)}\\
\hleq nK_1\left(c_8 \frac{M}{s \delta^\tau} + c_7\frac{M}{s}\right).
\end{eqnarray*}
The estimate
\[
 \left| \Delta W\right|_{\D(r_+, s_+)} = \max_{1\leq j \leq 2n} \left|
 \Delta W_j
 \right|_{\D(r_+, s_+)} \leq nK_1(c_7+c_8) \frac{M}{s\delta^\tau}\]
follows and the proof is finished.
\hfill$\Box$
\subsubsection*{Existence of the sequences}
Our intention is to formulate theorem
\ref{induktionslemma}
universally for the $k$-th step and to connect it with the Hamiltonian
(\ref{hami}).
To do that we have to find suitable sequences
$\left( r_k \right)$, $\left( \delta_k\right)$, $\left( s_k \right)$, and $\left( M_k\right)$. They shall allow it to use theorem
\ref{induktionslemma}
repeatedly with
\[ r=r_k, \; r_+ = r_{k+1},\; \delta = \delta_k,\; s= s_k,\;s_+ = s_{k+1},\mbox{ and }M=M_k.\]
At first we make sure that
$r_k$, $\delta_k$, and $s_k$ mesh correctly. We set
\begin{equation} \label{ansatzfolgen}
\delta_k := q^k \delta_0,\quad s_k := {\delta_k}^{\tau +1},\quad r_k := \frac{3}{4}r + 8 \delta_k\fa k\in \N_0,
\end{equation}
where $r$ is given in the assumptions of Theorem
\ref{hauptan},
$\delta_0 \in (0,1)$ is to be determined later, and
\begin{equation} \label{ansatzq}
q:=\frac{1}{4}.
\end{equation}
(\ref{ansatzfolgen}) yields immediately
\[ \delta_{k+1} = q^{k+1} \delta_0 = q \delta_k \mbox{ and } s_{k+1} =
   {\delta_{k+1}}^{\tau+1} = q^{\tau+1} s_k\fa k \in \N_0.\]
\begin{lemma} \label{folgenreigen}
The sequences 
$\left( r_k \right)_{k=0}^\infty$, $\left( \delta_k\right)_{k=0}^\infty$ and $\left( s_k \right)_{k=0}^\infty$
of
(\ref{ansatzfolgen}) and (\ref{ansatzq})
are decreasing and fulfill
\[ r_k > \frac{3}{4}r, \quad 0<\delta_k<\frac{r_k}{6},\quad 0<s_k\leq {\delta_k}^{\tau+1} \leq 1,\]\[0< r_{k+1} \leq r_k-6\delta_k,\quad 0<s_{k+1} \leq\frac{s_k}{8}\fa k\in \N_0.\]
\end{lemma}\noindent
{\bf Proof.} 
That the sequences decrease and that
$r_k > 3r/4$ for all $k\in  \N_0$ is clear.
We have
\[ \delta_k < \frac{8}{6} \delta_k < \frac{1}{6} \left( \frac{3}{4}r + 8
\delta_k\right) = \frac{r_k}{6}\fa k \in \N_0.\]
The definition of
$s_k$ and $\delta_k$ $(k\in \N_0)$
imply
$0 < s_k \leq {\delta_k}^{\tau +1} \leq 1$.
It is
$r_{k+1} = 3r/4 + 8 \delta_{k+1}$ and $r_k - 6\delta_k = 3r/4 + 2\delta_k$. 
Therefore
$r_{k+1} \leq r_k - 6 \delta_k$
holds if and only if
\[ 8 \delta_{k+1} \leq 2 \delta_k \gdw 4 q^{k+1} \delta_0 \leq q^k \delta_0
\gdw 4q \leq 1,\]
which is indeed true according to
(\ref{ansatzq}).
From
$\tau+1 \geq 2$
we infer
\[ s_{k+1} = \left( q^{k+1} \delta_0 \right)^{\tau +1} = q^{\tau +1} s_k \leq
q^2 s_k = \frac{s_k}{16} < \frac{s_k}{8}.\]
The Lemma is proved.
\hfill$\Box$\\\\
For the inductive lemma it is required to have sequences of functions
$\left(H_k\right)$, $\left( N_k\right)$, and $\left( R_k\right)$
which can be inserted for
$H$, $N$, and $R$, respectively, in the assumptions of theorem
\ref{induktionslemma}.
Let us suppose there are normal forms
$N_\ell$ 
defined on
$\D(r_\ell, s_\ell)$
$(0\leq \ell\leq k+1,\,k\in \N_0)$,
which meet 
(\ref{ind1Q+-Q})
and let us suppose 
$N_0$ fulfills something like
(\ref{nichtdegan}),
namely
\[ \left| N_{0yy}-C\right|_{\D(r_0,s_0)} \leq \frac{1}{4 |C^{-1}|}.\]
Then
\begin{eqnarray*}
\left| N_{k+1 \eta\eta}-C\right|_{\D(r_{k+1}, s_{k+1})} \hleq \sum_{\ell=0}^k \left| N_{\ell+1\eta\eta} - N_{\ell\eta\eta}
\right|_{\D(r_{\ell+1}, s_{\ell+1})} + \left| N_{0\eta\eta} - C \right|_{\D(r_0, s_0)}\\ \hleq
\sum_{\ell=0}^\infty c_{11} \frac{M_\ell}{{s_\ell}^2} + \frac{1}{4|C^{-1}|}
\end{eqnarray*}
is a consequence. Having
(\ref{ind1vorB}) in mind we therefore require
\begin{equation} \label{bedreihe}
\sum_{k=0}^\infty \frac{M_k}{{s_k}^2} \leq c_{17}, \quad c_{17} = \frac{1}{4 c_{11} |C^{-1}|}.
\end{equation}
From
(\ref{ind1vorC}) and (\ref{ind1R+})
the requirements
\begin{equation} \label{bedrest}
c_{15} \frac{{M_k}^2}{{s_k}^2} \leq M_{k+1} \mbox{ and } M_k \leq c_{18}
{s_k}^2\fa k \in \N_0,\; c_{18} = \frac{1}{16(c_7 + c_8)}
\end{equation}
follow. Observe, that
$c_{17}$ and $c_{18}$
depend on
$n$, $\tau$, $\gamma$, and $C$
only. 
In order to fulfill (\ref{bedrest})
we choose
\begin{equation} \label{ansatzM}
M_k := \frac{{s_k}^2}{c_{15}}t_k, \quad t_k := {t_0}^{\mu^k} \equiv {t_0}^{(\mu^k)}\fa k\in \N_0,
\end{equation}
with some
$t_0 \in (0,1)$,
and
\begin{equation} \label{ansatzmu}
\mu := \frac{3}{2}.
\end{equation}
(\ref{ansatzM}) gives promptly
\[ t_{k+1} = {t_0}^{\mu^{k+1}} = {t_0}^{\mu \cdot \mu^k} = {t_k}^\mu\fa k
\in \N_0.\]
\begin{remark}{\em \label{klnennermu2}
In formulas
(\ref{ansatzq}) and (\ref{ansatzmu})
any other value of
$q \in (0, 1/4]$
and
$\mu \in (1,2)$
would have done it equally well. 
\\
The parameter $\mu$ may be interpreted as the speed of convergence.
However, $\mu=2$ is not possible. This is due to the small divisors
(compare remarks
\ref{klnennerbem} (page
\pageref{klnennerbem}) and \ref{klnennermu} (page \pageref{klnennermu})).
}\end{remark}
\begin{lemma} 
The inequality
$c_{15}\cdot c_{18} \geq 1$ holds. \label{c15c18}
\end{lemma}\noindent
{\bf Proof.} 
We do the proof by tracing back the definition of
$c_{15}$.
At first, 
(\ref{defc12}) determines
\[ c_{12} = 1 + \frac{4 c_6 |C|}{\gamma} \geq 1.\]
Using 
$n\geq 2$, $|C|\,|C^{-1}|\geq |C\,C^{-1}|=1$,
and
(\ref{defc13}) we obtain
\[ c_{13} = 1+n \,c_{12} \left( 1 + 4 |C|\,|C^{-1}|\right) \geq 1 + 5n \geq 11.\]
Hence we have for
$c_{10}$
(see definition (\ref{defc10})) 
\[ c_{10} = 1 + n^2 |C| \left( 2c_7 +c_8\right) + c_{13} \geq 12.\]
The constant
$c_{15}$
was defined in
(\ref{defc15}),
this yields
\[ c_{15} = n(1+c_{10})(4c_7+c_8) \geq 26(c_7+c_8).\]
Now we calculate
\[ c_{15} \cdot c_{18} = \frac{c_{15}}{16(c_7+c_8)} \geq \frac{26}{16} \geq 1,\]
and the lemma is proven.
\hfill$\Box$
\begin{lemma} \label{reihet}
Let 
$m>1$ and $0<t<1$.
Then the estimate
\[ \sum_{k=0}^\infty t^{m^k} \leq \frac{t}{1-t^{m-1}}\]
holds.
\end{lemma} \noindent
{\bf Proof.}
Because of the equality
\[ \frac{t}{1-t^{m-1}} = t\sum_{k=0}^\infty \left( t^{m-1}\right)^k\]
it is sufficient to prove
\[ t \left( t^{m-1} \right)^k \geq t^{m^k} \gdw k(m-1) +1 \leq m^k =
(1+(m-1))^k \fa k\in \N_0.\]
This amounts to Bernoulli's inequality, which implies the assertion.
\hfill$\Box$
\begin{lemma} \label{reihetklein}
There exists a constant
$c_{19} = c_{19}(n, \tau, \gamma, C)>0$,
such that
the sequence
$\left( M_k \right)_{k=0}^\infty$ defined in
(\ref{ansatzM})
satisfies the conditions
(\ref{bedreihe}) and (\ref{bedrest})
for all
$t_0 \in (0, c_{19}]$. Moreover
\begin{equation} \label{bedreihet}
\sum_{k=0}^\infty \frac{M_k}{{s_k}^2}  \leq \frac{2}{c_{15}} t_0
\end{equation}
holds.
\end{lemma} \noindent
{\bf Proof.} 
By definition of the 
$t_k$
we see
$t_{k+1} = {t_k}^\mu$ $(k\in \N_0)$.
We require
$c_{19} \leq q^{(2\tau+2)/(2-\mu)}$,
than
$t_0 \leq q^{(2\tau+2)/(2-\mu)}$
follows.
The sequence of the 
$t_k$ decreases, so
$t_k \leq q^{(2\tau+2)/(2-\mu)}$ for all $k\in \N_0$.
This means
${t_k}^{2-\mu} \leq q^{2\tau+2}$ $(k\in \N_0)$.
Furthermore we have
\[ s_{k+1} = {\delta_{k+1}}^{\tau+1} = \left( q\cdot \delta_k\right)^{\tau+1}
= q^{\tau+1} s_k\fa k\in \N_0.\]
Hence we obtain
\[ c_{15} \frac{{M_k}^2}{{s_k}^2} = \frac{1}{c_{15}} {s_k}^2 {t_k}^2 =
\frac{1}{c_{15}} \frac{{s_{k+1}}^2}{q^{2\tau+2}} {t_k}^{2-\mu} {t_k}^\mu \leq
\frac{1}{c_{15}} {s_{k+1}}^2\, t_{k+1} = M_{k+1}\fa k\in \N_0.\]
This is the first inequality
(\ref{bedrest}). The second one
(\ref{bedrest})
is equivalent to
\[ t_k \leq c_{15}\cdot c_{18}\fa k\in \N_0.\]
This in turn is a consequence of lemma
\ref{c15c18}.
(\ref{ansatzq}) and (\ref{ansatzmu}) imply
$c_{19} < q = 1/4 = (1/2)^{1/(\mu-1)}$.
Hence
${t_0}^{\mu-1} \leq {c_{19}}^{\mu-1} \leq 1/2$,
and with
(\ref{ansatzM}) 
and lemma \ref{reihet} we get
\[ \sum_{k=0}^\infty \frac{M_k}{{s_k}^2} = \frac{1}{c_{15}} \sum_{k=0}^\infty
   {t_0}^{\mu^k} \leq \frac{1}{c_{15}} \frac{t_0}{1-{t_0}^{\mu-1}} \leq
   \frac{2}{c_{15}}t_0,\]
which is formula
(\ref{bedreihet}). Let us diminish 
$c_{19}$
by setting
\[ c_{19} := \min \left\{ q^{\frac{2\tau+2}{2-\mu}}, \frac{c_{15}c_{17}}{2}
   \right\},\]
then
$t_0 \leq c_{19} \leq c_{15} c_{17} /2$ and (\ref{bedreihet})
imply
(\ref{bedreihe}).
All assertions are shown.
\hfill$\Box$
\\\\
We define the constants in the assumptions of Theorem
\ref{hauptan}
as follows:
\begin{equation} \label{paragraph4}
c_1 := \min \left\{ c_{19},
\frac{c_{15}}{32n^2(c_7+c_8)\exp(c_{14} c_{17})} \right\}, \quad c_2 := \frac{1}{32^{2(\tau +1)} c_{15}}.
\end{equation}
To remind: So far we encountered the positive constants
$c_6$ to $c_{19}$.
The constants
$c_1$ and $c_2$
were defined right now, and the constants
$c_3$, $c_4$, and $c_5$ from the assertions of theorem
\ref{hauptan} will be determined later.
\begin{lemma} \label{defd0t0}
Let
$r$, $s$, $M$, and $\teta$
be the constants from theorem
\ref{hauptan} and set
\begin{equation} \label{defdotoformel}
\delta_0 := \frac{1}{32} s^{\frac{1}{\tau+1}}, \quad t_0 := \teta.
\end{equation}
Then
$r_0$, $s_0$ 
given by
(\ref{ansatzfolgen}),
and
$M_0$ 
from 
(\ref{ansatzM})
with
$k=0$, satisfy
\[ r_0 \leq r, \quad s_0 \leq s, \quad M_0 \geq M.\]
\end{lemma} \noindent
{\bf Proof.}
The fact
$s \leq r^{\tau+1}$
and the definition of
$\delta_0$ 
show
\[ r_0 = \frac{3}{4}r + 8 \delta_0 \leq \frac{3}{4}r + \frac{1}{4}
s^{\frac{1}{\tau+1}} \leq r.\]
Furthermore
\[ s_0 = {\delta_0}^{\tau+1} = \frac{s}{32^{\tau+1}} < s\]
follows.
For the claim $M_0 \geq M$
it is sufficient to prove
$M_0 \geq c_2 s^2 \teta$ because of
$M \leq c_2 s^2 \teta$. We have
\[ c_2s^2 \teta \leq \frac{1}{32^{2(\tau+1)}c_{15}}s^2 \teta = \frac{1}{c_{15}} \left(
\frac{s^{\frac{1}{\tau+1}}}{32} \right)^{2(\tau+1)}\cdot \teta = \frac{1}{c_{15}}
     {\delta_0}^{2(\tau+1)} \teta = \frac{1}{c_{15}} {s_0}^2t_0 = M_0,\]
which proves the lemma.
\hfill$\Box$
\begin{theorem} \label{induktionssatz}
\addcontentsline{toc}{subsubsection}{\hspace{1.325cm}Satz
     \ref{induktionssatz} (Induktionslemma)}
{\em \bf (inductive lemma)} 
Under the assumptions of theorem
\ref{hauptan}
and with the sequences
$\left( r_k \right)_{k=0}^\infty$,
$\left( \delta_k\right)_{k=0}^\infty$, $\left( s_k \right)_{k=0}^\infty$, and
$\left( M_k\right)_{k=0}^\infty$
fixed in
(\ref{ansatzfolgen}), (\ref{ansatzq}), (\ref{ansatzM}), (\ref{ansatzmu}), and
(\ref{defdotoformel})
the following holds for all
$k\in \N_0$:\\
There exist simple canonical transformations
\begin{equation} \label{gestaltzk}
Z_{k+1}:\D(r_{k+1}, s_{k+1}) \lto \D(r_k-5\delta_k, s_k/4), \quad Z_{k+1} -
\id \in \Pp_{2n}(r_{k+1}, s_{k+1}),
\end{equation}
such that the functions
\begin{equation} \label{secunum}
H_{k+1} := H_k \circ Z_{k+1} = H_0 \circ Z_1 \circ Z_2 \circ \ldots \circ
Z_{k+1} \quad \mbox{ with } H_0 := \left. H\right|_{\D(r_0, s_0)}
\end{equation}
are elements of the respective space
$\Pp(r_{k+1}, s_{k+1})$
and can be written as 
$H_{k+1} = N_{k+1} + R_{k+1}$ with $N_{k+1}$, $R_{k+1} \in \Pp(r_{k+1}, s_{k+1})$, and
\begin{equation} \label{gestalthk}
N_{k+1} (\xi, \eta) = a_{k+1} +
\skal{\omega}{\eta} + \Oo(|\eta|^2), \quad a_{k+1} \in \R.
\end{equation}
The following estimates hold for all
$k \in \N_0$:
\begin{equation} \label{indkZ}
\left| Z_{k+1, \zeta} \right|_{\D(r_{k+1}, s_{k+1})} \leq \exp \left( c_{14} \frac{M_k}{{s_k}^2}\right),
\end{equation}
\begin{equation} \label{indkZ-E}
\left| Z_{k+1,\zeta}-E_{2n}\right|_{\D(r_{k+1}, s_{k+1})} \leq c_{14}
\frac{M_k}{{s_k}^2} \exp \left( c_{14} \frac{M_k}{{s_k}^2} \right),
\end{equation}
\begin{equation} \label{indka+-a}
\left| a_{k+1}-a_k\right|\leq \widetilde{c}_9 \frac{M_k}{s_k},
\end{equation}
\begin{equation} \label{indkQ+-Q}
\left| N_{k+1\eta\eta} - N_{k\eta\eta} \right|_{\D(r_{k+1}, s_{k+1})} \leq
c_{11} \frac{M_k}{{s_k}^2} \quad \left(N_0 := \left. (H-R) \right|_{\D(r_0,
  s_0)}\right)\!,
\end{equation}
\begin{equation} \label{indkR+}
\left| R_{k+1} \right|_{\D(r_{k+1}, s_{k+1})} \leq c_{15}
\frac{{M_k}^2}{{s_k}^2}.
\end{equation}
Herein the constants
$\widetilde{c}_9$
and
$c_{11}$
are given by Theorem
\ref{satzlingl},
$c_{14}$ and $c_{15}$
by Theorem
\ref{induktionslemma}.
Moreover
$W_{k+1} := Z_1 \circ \ldots \circ Z_{k+1}$
fulfills
\begin{equation} \label{indkWzeta}
 \left| W_{k+1,\zeta} \right|_{\D(r_{k+1}, s_{k+1})} \leq \exp\left( c_{14}
 \sum_{\ell=0}^k \frac{M_\ell}{{s_\ell}^2} \right) \fa k \in \N_0,
\end{equation}
and
$\Delta W_{k+1} := W_{k+1} - W_k$ $(k \in \N)$, $\Delta W_1 := W_1 - \id$
satisfies
\begin{equation} \label{indkDeltaW}
\left| \Delta W_{k+1} \right|_{\D(r_{k+1}, s_{k+1})} \leq c_{20}
\frac{M_k}{s_k {\delta_k}^\tau} \fa k \in \N_0,
\end{equation} 
where
$c_{20} = c_{20} (n, \tau, \gamma, C)$ 
is a positive constant.
\end{theorem}\noindent
{\bf Proof.}
Clearly the proof is to be done by repeated use of theorem
\ref{induktionslemma}.
Lemma 
\ref{defd0t0}
shows
$\D(r_0, s_0) \subseteq \D(r,s)$.
So
$H_0$
can be defined as the restriction of the function
\[ H (x,y) = a + \langle\omega, y\rangle + \frac{1}{2} \langle y\cdot Q(x) , y\rangle + R(x,y)\]
of
(\ref{hami}) to
$\D(r_0, s_0)$.
We set
$a_0 := a$ and $R_0 := \left. R\right|_{\D(r_0, s_0)}$
with
$a$ and $R$
from
(\ref{hami}).
To summarize, we start the induction in accordance with
(\ref{secunum}) and (\ref{indkQ+-Q}) with\\
\begin{minipage}{10cm}
\[
H_0 = \left. H \right|_{\D(r_0, s_0)},\, R_0 = \left. R \right|_{\D(r_0,
   s_0)},\, a_0 = a \mbox{ and } N_0 = \left. (H-R)\right|_{\D(r_0, s_0)},\]\vspace{-7mm}
\[ \mbox{ where } H, \,R \mbox{ and } a \mbox{ are given by (\ref{hami}).} 
\]
\end{minipage}
\begin{minipage}{5.3cm}
\begin{equation}
\label{secu2}
\end{equation}
\end{minipage}\\\\
We check the assumptions of theorem
\ref{induktionslemma}.
The assumptions on the constants
$r$, $\delta$, $s$, $r_+$ and $s_+$
are fulfilled by lemma
\ref{folgenreigen}.
Apply the lemma for
$k=0$ and
\[ r = r_0, \; \delta = \delta_0, \; s = s_0, \; r_+
= r_1, \; s_+ = s_1.\]
In theorem 
\ref{induktionslemma}
we use
\[ H=H_0,\;N=N_0 = H_0-R_0, \,R=R_0\mbox{ and }M = M_0\]\vspace{-7mm}
\[ \mbox{with } H_0,\, N_0, \, R_0 \mbox{ from (\ref{secu2}) and $M_0$ from
(\ref{ansatzM}) for $k=0$}.\]
Then the function $N$ of 
(\ref{ind1vorA})
has the form
\[ N(x,y) = a_0 + \skal{\omega}{y} + \frac{1}{2} \skal{y\cdot Q(x)}{y}\fa
(x,y) \in \D(r_0, s_0)\]
because of
(\ref{hami}).
So
(\ref{nichtdegan}) implies
(\ref{ind1vorB}).
Lemma
\ref{defd0t0} and (\ref{Man})
show
\[ \left| R_0 \right|_{\D(r_0, s_0)} = \left| R \right|_{\D(r_0, s_0)} \leq
\left| R \right|_{\D(r,s)} =M \leq M_0.\]
Moreover by the inequality
(\ref{bedrest}),
which holds according to Lemma
\ref{reihetklein},
we have
\[ M_0 \leq c_{18} {s_0}^2 = \frac{1}{16(c_7+c_8)}{s_0}^2.\]
Hence assumption 
(\ref{ind1vorC})
is met and we may apply theorem
\ref{induktionslemma}.
It yields a transformation $Z$ and a function
$H_+$, as well as
$a_+$, $N_+$, and
$R_+$.
Now we set
\[ Z_1 := Z,\;H_1 := H_+ \in \Pp(r_1,s_1),\; a_1 := a_+ \in \R,\]\[N_1:= N_+
\in\Pp(r_1, s_1)\mbox{ and } R_1 := R_+ \in \Pp(r_1,s_1).\]
Then assertions
(\ref{gestaltzk})
to
(\ref{indkR+}) follow for $k=0$.
In case
$k=0$
(\ref{indkWzeta})
is equivalent to
(\ref{indkZ})
because of
$W_1 = Z_1$. Hence (\ref{indkWzeta}) holds.
To prove
(\ref{indkDeltaW})
for
$k=0$
we consider
$\Delta W_1 = Z_1 - \id = \id \circ Z_1 - \id$.
So let us put
$W=\id$
and $K_1 = 1$ in theorem
\ref{induktionslemma}, then we obtain with
(\ref{ind1DeltaW})
\[ \left| \Delta W_1\right|_{\D(r_1, s_1)} \leq n(c_7+c_8)\frac{M_0}{s_0
  {\delta_0}^\tau}.\]
We define
\begin{equation} \label{defc20}
c_{20} := n(c_7+c_8) \exp(c_{14}c_{17}),
\end{equation}
then (\ref{indkDeltaW}) holds for $k=0$. (The reason for the factor $\exp(c_{14} c_{17})$ will become clear at the end of the proof.)\\
Now suppose the inductive Lemma is true for all
$\ell$, $0 \leq \ell \leq k-1 \in \N_0$. We want to apply theorem
\ref{induktionslemma} with
\[ r= r_k, \; \delta = \delta_k, \;s= s_k, \; r_+
= r_{k+1}, \; s_+ = s_{k+1}. \]
Lemma
\ref{folgenreigen}
says that the assumptions on these constants are fulfilled. Next we have to put
\[ H=H_{k},\; a=a_{k},\;N=H_k-R_k, \mbox{ and }R=R_k.\]
By lemma
\ref{reihetklein},
formula (\ref{bedreihe})
holds, namely
\[ \sum_{k=0}^\infty \frac{M_k}{{s_k}^2} \leq \frac{1}{4c_{11}|C^{-1}|}.\]
Using
(\ref{indkQ+-Q})
up to
$k-1$
we get
\begin{eqnarray*}
\left| N_{k\eta\eta} -C\right|_{\D(r_k, s_k)} \hleq \sum_{\ell=0}^{k-1} \left| N_{\ell+1\eta\eta} - N_{\ell\eta\eta}
\right|_{\D(r_{\ell+1}, s_{\ell+1})} + \left| N_{0\eta\eta} -C\right|_{\D(r_0, s_0)} \\\hleq c_{11}
\frac{1}{4 c_{11} |C^{-1}|} + \frac{1}{4|C^{-1}|} = \frac{1}{2 |C^{-1}|}.
\end{eqnarray*}
So assumption
(\ref{ind1vorB})
is satisfied.
(\ref{bedrest}) holds because of lemma
\ref{reihetklein}, in particular we have
\[ c_{15} \frac{{M_{k-1}}^2}{{s_{k-1}}^2} \leq M_k \mbox{ and } M_k \leq
  \frac{1}{16(c_7+c_8)}{s_k}^2.\]
Hence
(\ref{indkR+})
for
$k-1$
shows
\[ \left| R_k \right|_{\D(r_k, s_k)} \leq c_{15} 
\frac{{M_{k-1}}^2}{{s_{k-1}}^2} \leq \frac{1}{16(c_7+c_8)} {s_k}^2,\]
this is assumption
(\ref{ind1vorC}).
Theorem \ref{induktionslemma} can be applied and yields a transformation $Z$ and a function
$H_+$, as well as
$a_+$, $N_+$, and $R_+$.
Now we set
\[ Z_{k+1} := Z,\;H_{k+1} := H_+ \in \Pp(r_{k+1},s_{k+1}),\; a_{k+1} := a_+ \in \R,\]\[N_{k+1}:= N_+
\in\Pp(r_{k+1}, s_{k+1})\mbox{ and } R_{k+1} := R_+ \in \Pp(r_{k+1},s_{k+1}).\]
Assertions 
(\ref{gestaltzk})
to
(\ref{indkR+})
follow for the index $k$. To prove
(\ref{indkWzeta}) we calculate
\[ W_{k+1,\zeta} = Z_{1\zeta}(Z_2 \circ \ldots \circ Z_{k+1})\cdot Z_{2\zeta} (Z_3
\circ \ldots \circ Z_{k+1})\cdot \ldots \cdot Z_{k+1,\zeta}.\]
Formula
(\ref{indkZ})
up to $k$ implies
\begin{eqnarray*}
\left| W_{k+1,\zeta} \right|_{\D(r_{k+1}, s_{k+1})} \hleq \left| Z_{1\zeta}
\right|_{\D(r_1, s_1)} \cdot \left| Z_{2\zeta} \right|_{\D(r_2, s_2)} \cdot
\ldots \cdot \left| Z_{k+1,\zeta} \right|_{\D(r_{k+1}, s_{k+1})} \\ \hleq
\prod_{\ell=0}^{k} \exp \left( c_{14} \frac{M_\ell}{{s_\ell}^2} \right) = \exp
\left( c_{14} \sum_{\ell=0}^k \frac{M_\ell}{{s_\ell}^2} \right),
\end{eqnarray*}
so
(\ref{indkWzeta})
is shown for the index $k$. Furthermore
(\ref{indkWzeta}) for
$k-1$
and (\ref{bedreihe}), which holds by Lemma
 \ref{reihetklein}, give the estimate
\[ \left| W_{k\zeta} \right|_{\D(r_k, s_k)} \leq \exp \left( c_{14}
\sum_{\ell=0}^{k-1} \frac{M_\ell}{{s_\ell}^2} \right) \leq \exp ( c_{14} c_{17}).\]
Therefore we can insert
 $K_1 = \exp(c_{14} c_{17})$ in formula
(\ref{ind1DeltaW})
and
(\ref{indkDeltaW})
follows for the index $k$. Altogether the inductive lemma is proved.
\hfill$\Box$
\setcounter{equation}{0}
\section{Convergence of the iterative process}
In this section we complete the proof of theorem
\ref{hauptan}. Henceforth we work with the {\em general assumption}:
\begin{quote}
{\em Let the assumptions of theorem \ref{hauptan} be fulfilled. Let the sequences $\left( r_k \right)_{k=0}^\infty$, $\left( \delta_k
  \right)_{k=0}^\infty$, $\left( s_k \right)_{k=0}^\infty$, and $\left( M_k
  \right)_{k=0}^\infty$ be defined according to (\ref{ansatzfolgen}), (\ref{ansatzq}),
  (\ref{ansatzM}), (\ref{ansatzmu}), and (\ref{defdotoformel}).}
\end{quote}
Especially lemmas 
\ref{folgenreigen}, \ref{reihetklein}, and
\ref{defd0t0},
and the inductive lemma
\ref{induktionssatz} hold under this general assumption. 
\subsubsection*{Convergence of the symplectic transformations}
\begin{theorem} \label{mazda}
The maps
\[ W_k = Z_1 \circ \ldots \circ Z_k \quad (k \in \N)\]
provided by theorem
\ref{induktionssatz}
are simple canonical transformations.
$W_k - \id \in \Pp_{2n}(r_k, s_k)$ holds.
\end{theorem}\noindent
{\bf Proof.}
The maps
$W_k$
are well-defined, for 
$Z_{k+1}$ 
lies in the domain of
$Z_k$
for all $k\in \N$ 
by
(\ref{gestaltzk}).
The $W_k$ are simple canonical transformations.
Moreover
$W_k - \id \in \Pp_{2n}(r_k, s_k)$ holds for all $k\in \N$.\hfill$\Box$\\\\
Simple canonical transformations are affine-linear in $\eta$, so they can always be defined for all
$\eta \in \C^n$.
More precisely, if
$W_k = (U_k, V_k)$
is defined on
$\D(r_k, s_k)$
by
\begin{equation} \label{sternwk} 
W_k(\xi,\eta) = (U_k(\xi), V_k(\xi,0) + \eta\cdot U_{k\xi}(\xi)^{-1})\fa
(\xi,\eta) \in \D(r_k,s_k),
\end{equation}
as it is seen in theorem
\ref{einftrafos1}, then there exists a simple canonical transformation
$\widetilde{W}_k$
defined on
$\Ss(r_k) \times \C^n$
with
$\left. \widetilde{W}_k \right|_{\D(r_k, s_k)} = W_k$.
The equation
\begin{equation} \label{sternwkschlange} 
\widetilde{W}_k(\xi,\eta) = (U_k(\xi), V_k(\xi,0) + \eta\cdot U_{k\xi}(\xi)^{-1})\fa
(\xi,\eta) \in \Ss(r_k) \times \C^n
\end{equation}
holds. Comparing
(\ref{sternwk})
and
(\ref{sternwkschlange})
we notice that
$\Wtilde_k(\stv,0) = W_k(\stv,0)$.
When we write
$\Wtilde_k = (\widetilde{U}_k, \widetilde{V}_k)$,
we have 
$\widetilde{U}_k = U_k$ and $\widetilde{V}_{k\eta} = V_{k\eta}$ too. 
We will use this in the sequel.
\begin{theorem} \label{teilkomp}
There exists a subsequence
$\left( \Wtilde_{k_\ell}\right)_{\ell=1}^\infty$
which converges uniformly on compact subsets of
$\Ss(3r/4) \times \C^n$
to a simple canonical transformation
$W_\infty$ with
$W_\infty - \id \in \Pp_{2n} (3r/4, s)$.
\end{theorem} \noindent
{\bf Proof.}
It is
$r_k > 3r/4$ for all $k \in \N$ by
(\ref{ansatzfolgen}).
Therefore all maps
$\Wtilde_k$ are defined for $\zeta \in \Ss(3r/4) \times \C^n$.
Looking at the assumptions of theorem
\ref{einftrafos2}
we calculate with
(\ref{indkDeltaW}),
$s_k \leq {\delta_k}^\tau$
(by lemma
\ref{folgenreigen}), and (\ref{bedreihe})
\begin{eqnarray*}
\sum_{k=0}^\infty \left| W_{k+1} - W_k \right|_{\Ss(3r/4) \times \{0\}}
\hleq \sum_{k=0}^\infty \left| \Delta W_{k+1} \right|_{\D(r_{k+1}, s_{k+1})} \\
  \hleq \sum_{k=0}^\infty c_{20} \frac{M_k}{s_k \delta_k^\tau} \leq c_{20}
  \sum_{k=0}^\infty \frac{M_k}{s_k^2} \leq c_{17} c_{20}.
\end{eqnarray*}
This means, that the functions
$\Wtilde_k(\stv,0)=W_k(\stv,0)$
converge uniformly on
$\Ss(3r/4)$, in particular they converge uniformly on compact subsets.
We use the row-sum norm, so 
(\ref{indkWzeta})
and
(\ref{bedreihe})
show
\[ \left| V_{k\eta} \right|_{\Ss(3r/4)} \leq \left| W_{k\zeta}
\right|_{\Ss(3r/4) \times \{0\}} \leq \exp (c_{14} c_{17})\fa k\in \N.\]
Hence the theorem of Montel
(see \cite{rothstein82}, theorem 1.6)
tells us that there exists a subsequence
$\left( V_{k_\ell,   \eta}\right)_{\ell=1}^\infty$
which converges uniformly on compact subsets of
$\Ss(3r/4)$.
Let us set
$\left( \Wtilde_{k_\ell} \right)_{\ell=1}^\infty$ and $\U = \Ss(3r/4)$
in the assumptions of theorem
\ref{einftrafos2}. Then the theorem may be applied and predicates, that the
sequence
$\left( \Wtilde_{k_\ell} \right)_{\ell=1}^\infty$
converges uniformly on compact subsets of
$\Ss(3r/4) \times \C^n$
against a simple canonical transformation
\begin{equation*}
W_\infty = \left( U_\infty, V_\infty \right) : \Ss(3r/4) \times \C^n \lto
 \C^{2n}.
\end{equation*}
The functions
$\Wtilde_{k_\ell}$
map real vectors to real values 
and the
$\Wtilde_{k_\ell} -  \id$
are
$2\pi$-periodic
by
(\ref{sternwk}), (\ref{sternwkschlange}),
and theorem
 \ref{mazda}.
Therefore we obtain
$W_\infty-\id \in \Pp_{2n}(3r/4, s)$
and the proof is finished.
\hfill$\Box$
\begin{theorem} \label{www}
The function 
$W_\infty$
of theorem
\ref{teilkomp}
fulfills
\begin{equation} \label{Wunabb}
W_\infty (\zeta) \in \D(r,s)\fa \zeta \in \D(r/2, 5s/8).
\end{equation}
The restriction
\[ \fbox{$ \displaystyle W = (U, V) := \left. W_\infty \right|_{\D(r/2,
    s/2)}$}\]
is a simple canonical transformation with
\[ W: \D(r/2, s/2) \lto \D(r,s), \quad W-\id \in \Pp_{2n}(r/2, s/2).\]
There exists a positive constant
$c_3$,
which depends on
$n$, $\tau$, $\gamma$, and $C$
only, such that
\[ \left| W_\zeta - E_{2n} \right|_{\D(r/2, s/2)} \leq c_3 \teta.\]
\end{theorem}\noindent
{\bf Proof.}
The definition of $W$ and theorem
\ref{teilkomp} show that
$W-\id \in \Pp_{2n}(r/2, s/2)$
and that 
$W$ is a simple canonical transformation.
By the definition in theorem
\ref{induktionssatz}
we have
\[ W_k = Z_1 \circ \ldots \circ Z_k \quad (k \in \N).\]
Let us write
$W_k = (U_k, V_k)$.
The functions
$Z_k = (X_k, Y_k)$
are simple canonical transformations, so 
\[ U_k = U_k(\xi) = X_1 \circ \ldots \circ X_k(\xi).\]
In particular the functions
$U_k$ map to
$\Ss(r_0 - 5\delta_0)$ by
(\ref{gestaltzk}).
The function $U$
is the limit of a subsequence of the
$U_k$. Hence $U$ is defined on
$\Ss(r/2)$
and maps to
$\Ss(r_0 - 4\delta_0)$.
Because of
lemma \ref{defd0t0} $r_0 \leq r$, so $\Ss(r_0 - 4\delta_0) \subseteq \Ss(r)$, and consequently
\[ U: \Ss(r/2) \lto \Ss(r).\]
By definition of
$W$
we have
$U = \left. U_\infty \right|_{\Ss(r/2)}$. This implies
\[U_\infty(\xi) \in \Ss(r)\fa \xi \in \Ss(r/2).\]
Next (notice
(\ref{Wunabb})) we have to prove
\[ \left| V_\infty(\xi,\eta) \right| < s \fa (\xi,\eta) \in \D(r/2, 5s/8). \]
To that end we observe for
$(\xi,\eta) \in \D(3r/4, 5s/8)$ 
\begin{equation} \label{dokov}
V_\infty(\xi,\eta) = V_\infty(\xi,0) + \eta\, U_{\infty\xi}(\xi)^{-1} = V_\infty(\xi,0) +\eta + \eta \left(
U_{\infty\xi}(\xi)^{-1} - E_n \right).
\end{equation}
We consider
$V(\stv,0)$.
Each
$W_k$ $(k \in \N)$
maps
$(\xi,0) \in \Ss(r_k) \times \{0\}$
to
$\D(r_0 - 5\delta_0, s_0/4)$,
for this is true for
$Z_1$.
Therefore
$\left| V_k(\stv,0)\right|_{\Ss(r_k)} < s_0/4$
holds for all
$k \in \N$.
This implies
$\left| V(\stv,0)\right|_{\Ss(r/2)} \leq s_0/4$,
and with
$s_0 \leq s$
(by lemma
\ref{defd0t0})
we obtain
\begin{equation} \label{dopfuv}
\left| V(\xi,0) \right| \leq \frac{s}{4}\fa \xi \in \Ss(r/2).
\end{equation}
We need an estimate for
$U_\xi^{-1} - E_n$.
It can be found with lemma
\ref{invmatrix}.
Thereto we search for an inequality for
$U_\xi - E_n$.
We have for all
$k\in \N$
and all
$\zeta \in \D(r_k, s_k)$ 
\begin{equation} \label{bmeise2}
W_k(\zeta) - \zeta = \Delta W_1(\zeta) + \ldots + \Delta W_k (\zeta).
\end{equation}
(\ref{indkDeltaW})
and Cauchy's estimate show for
$k \in \N_0$
\[ \left| \Delta W_{k+1,\xi} \right|_{\D(r_{k+1}-\delta_k, s_{k+1})} \leq c_{20}
\frac{M_k}{s_k {\delta_k}^\tau \cdot \delta_k} \leq c_{20} \frac{M_k}{{s_k}^2}.\]
By
(\ref{ansatzfolgen}) and (\ref{ansatzq})
we see
\[ r_{k+1} - \delta_k = \frac{3r}{4} + 8 \delta_{k+1} - \delta_k =
\frac{3r}{4} + 8q \delta_k - \delta_k = \frac{3r}{4} + \delta_k >
\frac{3r}{4}\fa k \in \N_0.\]
So
(\ref{bedreihet}) and (\ref{bmeise2})
yield the estimate
\begin{equation} \label{ameise2} 
\left| W_{k\xi} - \left( \!\!\begin{array}{c} E_n \\ 0 \end{array} \!\!\right) \right|_{\Ss(3r/4) \times \{0\}} \leq
\sum_{\ell=0}^\infty \left| \Delta W_{\ell+1,\xi} \right|_{\D(r_{\ell+1}-\delta_\ell,
  s_{\ell+1})} \leq \frac{2 c_{20}}{c_{15}} t_0.
\end{equation}
Let us write
$\Delta W_k = (\Delta U_k, \Delta V_k)$.
Then in particular
\[ \left| U_{k\xi} - E_n \right|_{\Ss(3r/4)} \leq
\sum_{\ell=0}^\infty \left| \Delta U_{\ell+1,\xi} \right|_{\Ss(r_{\ell+1}-\delta_\ell)} \leq \frac{2 c_{20}}{c_{15}} t_0 \]
follows (note that we use the row-sum norm).
When we have a look at
(\ref{paragraph4}) and (\ref{defc20}), we see
\[ c_1 \leq \frac{c_{15}}{32 n c_{20}}.\]
It is $t_0 = \teta$ by (\ref{defdotoformel}) and $\teta \leq c_1$ by assumption of theorem \ref{hauptan}, so
\begin{equation} \label{ameise1}
\left| U_{k\xi} - E_n \right|_{\Ss(3r/4)} \leq \frac{2 c_{20}}{c_{15}} \teta
\leq \frac{1}{16n} \leq \frac{1}{16} \fa k \in \N.
\end{equation}
Now we can apply lemma
\ref{invmatrix}.
Therein we have to put
$S= E_n$, $P= U_{k\xi}(\xi)$ $(\xi \in \Ss(3r/4))$
and
$h = 2c_{20} \teta / c_{15}$.
The lemma says that
$U_{k\xi}(\xi)^{-1}$
satisfies the estimate
\begin{equation} \label{bmeise3}
\left| U_{k\xi}(\xi)^{-1} - E_n \right| \leq \frac{2c_{20}}{c_{15}} \,\teta\,
\frac{1}{1-\frac{1}{16}} = \frac{16}{15} \frac{2 c_{20}}{c_{15}}\,\teta \leq
\frac{1}{15n}\fa \xi \in \Ss(3r/4).
\end{equation}
This implies
\[ \left| U_{\infty \xi}^{-1} -E_n \right|_{\Ss(3r/4)} \leq \frac{1}{15n} \quad
\mbox{and} \quad \left| U_{\xi}^{-1} -E_n \right|_{\Ss(r/2)} \leq \frac{1}{15n},
\]
which in turn together with
(\ref{dokov})
and
(\ref{dopfuv})
leads to
\[ | V_\infty(\xi,\eta)| < \frac{s}{4} + \frac{5\,s}{8} + \frac{5\,s}{8} n
\frac{1}{15n} = \frac{30+75+5}{120}\,s <s\fa (\xi,\eta) \in \D(r/2, 5s/8).\]
We obtain
\[ W_\infty(\xi,\eta) \in \D(r,s)\fa(\xi,\eta) \in \D(r/2, 5s/8),\]
as well as
\[ W:\D(r/2, s/2) \lto \D(r,s).\]
In order to find an inequality for
$\left| W_\zeta - E_{2n}\right|$
we observe
\[ W_\zeta - E_{2n} = \left( \begin{array}{cc} U_\xi-E_n & 0 \\ V_\xi &
  \left(U_\xi^{-1}\right)^\tp -E_n \end{array} \right).\]
(\ref{ameise1})
gives
\begin{equation} \label{ameise5}
\left| U_\xi - E_n \right|_{\Ss(r/2)} \leq \frac{2 c_{20}}{c_{15}} \,\teta,
\end{equation}
and
(\ref{bmeise3})
shows
\begin{eqnarray} \nonumber
\left| \left( U_\xi^{-1}\right)^\tp - E_n \right|_{\Ss(r/2)} \heq \left| \left(
  U_\xi^{-1} - E_n\right)^\tp \right|_{\Ss(r/2)} \\ \hleq n \left| U_\xi^{-1}
  -E_n\right|_{\Ss(r/2)} \leq n\,\frac{16}{15} \frac{2 c_{20}}{c_{15}}\, \teta
  < 3n\frac{c_{20}}{c_{15}}\,\teta. \label{ameise4}
\end{eqnarray}
Let's turn to
$V_\xi$.
By definition
$V = \left. V_\infty \right|_{\D(r/2,   s/2)}$ holds, and
\[ V_\infty(\xi,\eta) = V_\infty(\xi,0) + \left( V_\infty (\xi,\eta) -
  V_\infty(\xi,0)\right)\fa (\xi,\eta) \in \D(3r/4, s).\]
Hence with
(\ref{dokov})
we obtain
\begin{eqnarray} \nonumber
V_{\infty\xi}(\xi,\eta) \heq V_{\infty\xi}(\xi,0) + \frac{\partial}{\partial\xi} \left( V_\infty (\xi,\eta) -
  V_\infty(\xi,0) - \eta \right)\\
\heq V_{\infty \xi}(\xi,0) + \frac{\partial}{\partial\xi} \left( \eta \left( U_{\infty
  \xi}(\xi)^{-1} - E_n\right)\right) \;\;\forall\;\; (\xi,\eta) \in \D(3r/4,s). \label{ameise3}
\end{eqnarray}
From
(\ref{ameise2})
it follows with
$t_0 = \teta$
\[ \left| V_{k\xi}(\stv,0) \right|_{\Ss(3r/4)} \leq 2 \frac{c_{20}}{c_{15}}\,
\teta\fa k \in \N.\]
This inequality holds for the limit
$V_\infty$ as well and hence for
$V$ giving
\[ \left| V_\xi(\stv,0) \right|_{\Ss(r/2)} \leq 2 \frac{c_{20}}{c_{15}}\,\teta.\]
To get the second summand of
(\ref{ameise3}) under control we define
\[ u: \D(3r/4, s) \lto \C^n,\quad (\xi,\eta) \mapsto u(\xi,\eta) = \eta\left(
U_{\infty\xi}(\xi)^{-1} -E_n\right).\]
From
(\ref{bmeise3})
we see
\[ \left| U_{\infty\xi}^{-1} - E_n \right|_{\Ss(3r/4)} \leq \frac{32}{15}
\frac{c_{20}}{c_{15}} \, \teta \leq 3 \frac{c_{20}}{c_{15}} \, \teta,\]
which implies
\[ \left|u\right|_{\D(3r/4, s)} \leq s n\left| U_{\infty \xi}^{-1} -E_n
\right|_{\Ss(3r/4)} \leq 3 n\frac{c_{20}}{c_{15}}\,s\teta.\]
Hence Cauchy's estimate and
$s \leq r^{\tau+1} \leq r$
show
\[ \left| u_\xi \right|_{\D(r/2, s)} \leq 3 n\frac{c_{20}}{c_{15}}
\frac{4s}{r}\, \teta \leq 12 n\frac{c_{20}}{c_{15}}\,\teta.\]
Therefore we can conclude with
(\ref{ameise3}) that
\[ \left| V_\xi \right|_{\D(r/2, s/2)} \leq \left|
V_\xi(\stv,0)\right|_{\Ss(r/2)} + \left| u_\xi \right|_{\D(r/2, s)} \leq 2
\frac{c_{20}}{c_{15}}\, \teta + 12 n\frac{c_{20}}{c_{15}}\,\teta \leq
13 n \frac{c_{20}}{c_{15}}\,\teta.\]
For matrices we use the row-sum norm, so this estimate,
(\ref{ameise5}), and (\ref{ameise4})
yield
\begin{eqnarray*}
\lefteqn{\left| W_\zeta - E_{2n} \right|_{\D(r/2, s/2)} \leq \left| \left( \begin{array}{cc} U_\xi-E_n & 0 \\ V_\xi &
  \left(U_\xi^{-1}\right)^\tp -E_n \end{array} \right) \right|_{\D(r/2,
  s/2)}}\\
\hleq \max\left\{ \left| U_\xi - E_n \right|_{\D(r/2, s/2)},\, \left| V_\xi
  \right|_{\D(r/2,s/2)} + \left| \left( U_\xi^{-1}\right)^\tp - E_n
  \right|_{\D(r/2, s/2)} \right\}\\
 \\ \hleq (3n+13n) \frac{c_{20}}{c_{15}}\,\teta=c_3\teta,
\end{eqnarray*}
where
\[ c_3 = 16 n \frac{c_{20}}{c_{15}}\]
is a positive constant. The theorem is proved.
\hfill$\Box$
\subsubsection*{Proof of the properties of the transformed Hamiltonian}
\begin{theorem} \label{veRschwindibus}
The functions
$R_k$ $(k \in \N)$
provided by theorem
\ref{induktionssatz}
fulfill
\[ \left| R_k \right|_{\Ss(r/2)\times \{0\}} \lto 0, \; \left| R_{k\eta}
\right|_{\Ss(r/2)\times \{0\}} \lto 0 \mbox{ and } \left| R_{k\eta\eta}
\right|_{\Ss(r/2)\times \{0\}} \lto 0 \quad (k\to \infty).\]
\end{theorem}\noindent
{\bf Proof.}
The estimates 
(\ref{bedrest}) and (\ref{indkR+})
imply
\[ \left| R_k \right|_{\D(r_k, s_k)} \leq c_{15}
\frac{{M_{k-1}}^2}{{s_{k-1}}^2} \leq M_k\fa k\in \N.\]
From this we conclude with Cauchy's estimates
\[ \left| R_{k\eta} \right|_{\D(r_k, s_k/2)} \leq \frac{2M_k}{s_k}, \quad
\left| R_{k\eta\eta}\right|_{\D(r_k, s_k/4)} \leq \frac{8 M_k}{{s_k}^2}\fa k
\in\N.\]
The series
$\sum_{k=0}^\infty M_k/{s_k}^2$ is convergent, hence the sequences
$\left( M_k \right)_{k=0}^\infty$, $\left( 2M_k/s_k\right)_{k=0}^\infty$, and
$\left( 8 M_k/{s_k}^2 \right)_{k=0}^\infty$ tend to zero. This proves the theorem.
\hfill$\Box$
\begin{theorem} \label{tayloreisatz}
Let $H$ be the function of theorem
\ref{hauptan}. Then there exists a number
$a_+ \in \R$
and a function
$Q_+ \in \Pp_{n\times n} (r/2)$,
such that the Taylor expansion of
$H\circ W:\D(r/2, s/2) \lto \C$
is given by
\begin{equation} \label{taylorei}
H\circ W(\xi,\eta) = a_+ + \skal{\omega}{\eta} +\frac{1}{2} \skal{\eta \cdot
  Q_+(\xi)}{\eta} + \Oo(|\eta|^3).
\end{equation}
\end{theorem}\noindent
{\bf Proof.}
By theorem
\ref{induktionssatz}
we have
\begin{equation} \label{ueberfluss}
H_k = H \circ W_k = N_k + R_k  \fa k \in \N.
\end{equation}
So 
\[ H\circ W (\xi, 0) = \lim_{\ell\to\infty} H\circ W_{k_\ell} (\xi, 0) =
\lim_{\ell\to\infty} \left( N_{k_\ell}(\xi,0) + R_{k_\ell}(\xi,0)\right)\]
holds for all
$\xi \in \Ss(r/2)$.
The sequence
$R_{k_\ell}(\xi, 0)$
has the limit zero as we have seen in the theorem above. 
The sequence
$N_{k_\ell}(\xi, 0) = a_{k_\ell}$
is convergent because of
(\ref{indka+-a}),
we call its limit
\[ \fbox{$\displaystyle a_+ := \lim_{\ell\to\infty} a_{k_\ell}$}\]
The number
$a_+$ is a limit of real numbers, so it is a real number as well. We have
\[ H\circ W(\xi,0) = a_+\fa \xi \in \Ss(r/2).\]
Moreover we obtain for all
$\xi \in \Ss(r/2)$ by
(\ref{ueberfluss})
\begin{eqnarray*}
(H\circ W)_\eta(\xi,0) \heq H_z(W(\xi,0))\cdot W_\eta(\xi,0) = \lim_{\ell\to
    \infty} H_z(W_{k_\ell}(\xi,0))W_{{k_\ell},\eta} (\xi,0)\\
\heq \lim_{\ell\to\infty} \left( H\circ W_{k_\ell}\right)_\eta(\xi,0) =
    \lim_{\ell\to\infty} \left( N_{k_\ell,\eta}(\xi,0) + R_{k_\ell,\eta}(\xi,0) \right)= \omega.
\end{eqnarray*}
Now, the derivatives
$N_{k_\ell,\eta\eta}$
converge on
$\Ss(r/2) \times \{0\}$
by
(\ref{indkQ+-Q})
and we obtain a limit
\[ \fbox{$\dis Q_+(\xi) := \lim_{\ell\to\infty} N_{k_\ell,\eta\eta}(\xi,0)\fa \xi \in
\Ss(r/2)$}\]
This convergence is uniformly on 
$\Ss(r/2)$ and all functions
$N_{k_\ell,\eta\eta}(\stv,0)$ are elements of $\Pp_{n\times n}(r/2)$, so
$Q_+ \in \Pp_{n\times n}(r/2)$. Theorem 
\ref{teilkomp}
implies
\[  W_{k_\ell} (\stv,0) \longrightarrow W(\stv,0) \quad \mbox{uniformly on compact subsets of } \Ss(r/2).\]
Hence we conclude using the continuity of $W(\stv,0)$ and $H$
\[ H\circ W_{k_\ell}(\stv,0) \longrightarrow H\circ W(\stv,0) \quad
\mbox{uniformly on compact subsets of } \Ss(r/2).\]
Hence
(\ref{ueberfluss})
and the theorem of 
Weierstrass (see \cite{dieudonne60}, (9.12.1))
show 
for all $\xi \in \Ss(r/2)$
\begin{eqnarray*}
( H\circ W)_{\eta\eta} (\xi, 0) \heq \lim_{\ell\to\infty} \left( H\circ
W_{k_\ell}\right)_{\eta\eta}(\xi,0) =\lim_{\ell\to\infty}
\left( N_{k_\ell,\eta\eta}(\xi,0) + R_{k_\ell,\eta\eta}(\xi,0) \right)\\
\heq Q_+(\xi),
\end{eqnarray*}
which proves 
(\ref{taylorei}).
\hfill$\Box$
\begin{theorem} \label{hessesatz}
There exists a constant
$c_4 = c_4(n,\tau,\gamma, C)>0$, such that the function
$Q_+$ meets inequality
(\ref{hesse}), namely
\[ \left| Q_+  - Q\right|_{\Ss(r/2)} \leq c_4 \teta. \]
\end{theorem} \noindent
{\bf Proof.}
With
(\ref{indkQ+-Q}), (\ref{bedreihet}), $t_0 = \teta$, and
the fact that
$N_{0,\eta\eta}(\xi,0) = Q(\xi)$ holds for all $\xi \in \Ss(r/2)$
by definition of
$N_0$
in theorem
\ref{induktionssatz}, we conclude that
\[ \left| Q_+ - Q \right|_{\Ss(r/2)} \leq \sum_{k=0}^\infty c_{11}
\frac{M_k}{{s_k}^2} \leq \frac{2c_{11}}{c_{15}}\,\teta.\]
So, with the definition
\[ c_4 = \frac{2c_{11}}{c_{15}},\]
(\ref{hesse}) is shown.
\hfill$\Box$\begin{theorem} \label{restsatz}
There exists a number
$c_5 = 512/25 >0$, such that the function
\begin{equation} \label{defrstern}
R^{\ast}(\xi,\eta) :=(H\circ W)(\xi,\eta) - \left( a_+ + \skal{\omega}{\eta} +\frac{1}{2} \skal{\eta \cdot
  Q_+(\xi)}{\eta}\right),
\end{equation}
defined for all
$(\xi,\eta) \in \D(r/2, s/2)$, fulfills estimate
(\ref{tayl3}).
\end{theorem}\noindent
{\bf Proof.}
At first we observe that
$H\circ W_\infty(\xi,\eta)$
can be defined 
for all $(\xi,\eta) \in \D(r/2, 5s/8)$
by theorem \ref{www}. This gives an analytic continuation of 
$H\circ W$
to the domain
$\D(r/2, 5s/8)$. 
We call it
$H^{\ast\ast}$.
Therefore we can
enlarge definition
(\ref{defrstern})
to
$\D(r/2, 5s/8)$ and obtain an analytic continuation
$R^{\ast\ast}$ of
$R^\ast$.
Clearly 
(\ref{tayl3})
is equivalent to
\[ |R^{\ast\ast}(\xi, \eta)| \leq c_5 M \frac{|\eta|^3}{s^3} \quad
\mbox{for all } (\xi, \eta) \in \D(r/2, s/2), \]
which will be shown in the following. The derivatives with respect to
$\eta$
of
$H\circ W$
and
$H^{\ast\ast}$ coincide for all $(\xi,0) \in \Ss(r/2)\times \{0\}$.
So 
$R^{\ast\ast}(\xi,\eta) = \Oo(|\eta|^3)$ holds
by theorem
\ref{tayloreisatz}.
Moreover 
$R^{\ast\ast}$
is an analytic function. We fix an arbitrary
$\xi \in \Ss(r/2)$,
set
$N := H-R$, and consider
\[ H^{\ast\ast}(\xi,\eta) = H \circ W_\infty(\xi,\eta) = N\circ W_\infty(\xi,\eta) +
R\circ W_\infty(\xi,\eta) \quad (|\eta| < 5s/8).\]
Well, 
$W_\infty(\xi,\eta)$
is a polynomial of degree one in
$\eta$ and $N$
is, by
(\ref{hami}), a polynomial of degree two in
$\eta$.
Therefore
$N \circ W_\infty(\xi,\eta)$
has degree two in 
$\eta$
and the terms of order three and higher in
$\eta$
of
$H^{\ast\ast}(\xi,\stv)$
and
$R\circ W_\infty(\xi,\stv)$
coincide. Hence the same holds for
$R^{\ast\ast}(\xi,\stv)$
and
$R\circ W_\infty(\xi, \stv)$.
So we can apply lemma
\ref{ttaysatz},
in which the function 
$\eta \mapsto R\circ W_\infty(\xi,\eta)$
is bounded by
$M$
for $|\eta| < 5s/8$
because of
(\ref{Man})
and theorem
\ref{www}.
Putting
\[ \sigma = \frac{5s}{8},\;f=R\circ W_\infty(\xi,\stv) \mbox{ and } \eps =
\frac{4}{5}\]
in lemma
\ref{ttaysatz}, we obtain
\[ \left| R^{\ast\ast} (\xi,\eta)\right| \leq 5M\frac{|\eta|^3}{(5s/8)^3} =
\frac{512}{25} M \frac{|\eta|^3}{s^3}\fa |\eta| < \frac{4}{5} \frac{5s}{8} =
\frac{s}{2}.\]
Now,
$\xi \in \Ss(r/2)$ was arbitrary, so (\ref{tayl3}) holds with
\[ c_5 = \frac{512}{25}\]
and the theorem is proved.
\hfill$\Box$\\\\
Altogether theorems
\ref{www}, \ref{tayloreisatz},
\ref{hessesatz}, and \ref{restsatz}
prove theorem
\ref{hauptan}.
\begin{appendix}
\renewcommand{\theequation}{\Alph{section}.\arabic{equation}}
\renewcommand{\thetheorem}{\Alph{section}.\arabic{theorem}}
\setcounter{equation}{0}
\section{Appendix}
\subsection{A lemma on non-singular matrices}
\begin{lemma} \label{invmatrix}
Let
$S \in \C^{n\times n}$
be an invertible matrix. Then each matrix
$P\in \C^{n\times n}$ with
\[ |P-S| \leq h\cdot \frac{1}{|S^{-1}|}, \qquad 0 < h <1,\]
is invertible as well. The inverse of
$P$ fulfills
\[ |P^{-1}| \leq \frac{|S^{-1}|}{1-h} \quad \mbox{ and } \quad |P^{-1} - S^{-1}| \leq
\frac{h |S^{-1}|}{1-h}.\]
\end{lemma} \noindent
{\bf Proof.}
We set
$H := E_n - S^{-1}P$. The assumption leads to the estimate
\[ |H| = |E_n - S^{-1}P| \leq |S^{-1}|\,|S-P|\leq h < 1.\]
Therefore the Neumann series
\[ \sum_{k=0}^\infty H^k = (E_n -H)^{-1} = (S^{-1}P)^{-1}\]
converges, in particular
$S^{-1}P$ is non-singular. Hence this is also true for
$P = S\cdot S^{-1}P$. For
$P^{-1} = (S^{-1}P)^{-1}S^{-1}$
we find the estimate
\[ |P^{-1} | \leq |S^{-1}| \sum_{k=0}^\infty |H|^k \leq \frac{|S^{-1}|}{1-h}.\]
For
$P^{-1} - S^{-1} = (P^{-1}S-E_n)S^{-1}$
we calculate
\[ P^{-1} S - E_n = \left( \sum_{k=0}^\infty H^k \right) - E_n = \sum_{k=1}^\infty H^k\]
to see
\[ |P^{-1} - S^{-1} | \leq |S^{-1}| \sum_{k=1}^\infty |H|^k \leq \frac{h
  |S^{-1}|}{1-h},\]
as was to be shown.
\hfill $\Box$
\subsection{Estimates for analytic maps} \label{anc}
\begin{definition}{\em \label{defum}
Let
$z \in \C^n$ and $s>0$. We set
\[ \B(s;z) := \left\{ y \in \C^n\,|\, |y-z| < s\right\}.\]
}\end{definition}
The following lemma is Cauchy's estimate for analytic functions of several variables.
\begin{lemma} \label{cauchysabs}
Let
$M>0$ and $f:\B(s;0) \subseteq \C^n \to \C^m$
be an analytic function with
\[ |f|_{\B(s;0)} \leq M.\]
The the Jacobian of $f$ satisfies the estimate
\[ |f_x|_{\B(s-\eps;0)} \leq \frac{M}{\eps}
\mbox{ for all }
0<\eps < s.\]
\end{lemma}\noindent
{\bf Proof.} We fix an arbitrary $x_0 \in \B(s-\eps;0)$. Then (\ref{opnorm}) shows
\[ |f_x(x_0)| = \max_{|y|=1} |y f_{\,x}^{\rm T}(x_0)| = \max_{1\leq k \leq m} \max_{|y|=1}|\skal{f_{kx}(x_0)}{y}|,\]
where
$f_k$
denotes the $k$-th coordinate function of $f$. We give us arbitrary
$k \in \{1, \ldots, m\}$
and
$y \in \C^n$
with
$|y|=1$
and consider the auxiliary function
\[ g: \B(\eps;0) \subseteq \C \longrightarrow \C, \quad t \mapsto f_k(x_0 + ty).\] 
We obtain
\[ g_t(t) = \skal{f_{kx}(x_0 + ty)}{y} \quad \Rightarrow \quad g_t(0) = \skal{f_{kx}(x_0)}{y},\]
and Cauchy's estimate in one dimension says
\[ |\skal{f_{kx}(x_0)}{y}| = | g_t(0)| \leq \frac{M}{\eps},\]
which finishes the proof.
\hfill $\Box$\\\\
We need an estimate for the remainder of order three relating to the Taylor expansion of an analytic function. At first we prove it in dimension one.
\begin{lemma} \label{taylor3dim1}
Let
$\sigma >0$
and
$ g: \B(\sigma;0) \subseteq \C\to \C$, $z\mapsto g(z)$
be an analytic function bounded by a constant
$M>0$. Then the remainder
\[ h^{(g)}(z) := \sum_{k=3}^\infty \frac{1}{k!} \frac{\partial^k g}{\partial z^k}(0) \,z^k \quad (|z|<\sigma)\]
satisfies for all
$\eps\in (0,1)$
the estimate
\[|h^{(g)}(z)| \leq \frac{M}{1-\eps}\frac{|z|^3}{\sigma^3}\fa |z| \leq \eps\sigma.\]
\end{lemma}\noindent
{\bf Proof.}
By Cauchy's formula we have for
$0<\tsigma<\sigma$ 
\[ \left| \frac{\partial^k g}{\partial z^k}(0) \right| = \left| \frac{k!}{2\pi i} \oint_{|z|=\tsigma}
\frac{g(z)}{z^{k+1}}\,dz\right| \leq \frac{M k!}{\tsigma^k}.\]
The limit
$\tsigma \to \sigma$
yields
\[ \left| \frac{\partial^k g}{\partial z^k}(0)\right| \leq \frac{Mk!}{\sigma^k}.\]
Hence we get for the remainder, in case
$|z| \leq \eps \sigma$,
\begin{eqnarray*}
\left|h^{(g)}(z)\right| \hleq \sum_{k=3}^\infty \frac{1}{k!} \left|
\frac{\partial^k g}{\partial z^k}(0) \right| |z|^k \leq \sum_{k=3}^\infty
\frac{1}{k!} \frac{Mk!}{\sigma^k} |z|^k = M \sum_{k=3}^\infty \left(
\frac{|z|}{\sigma}\right)^k \\ \heq M \left( \frac{|z|}{\sigma}\right)^3
\sum_{k=0}^\infty \left( \frac{|z|}{\sigma}\right)^k \leq M \left(
\frac{|z|}{\sigma}\right)^3 \sum_{k=0}^\infty \eps^k = \frac{M}{1-\eps}
\frac{|z|^3}{\sigma^3},
\end{eqnarray*}
as was to be shown.
\hfill$\Box$
\begin{lemma} \label{ttaysatz}
Let
$\sigma >0$ and $f:\B(\sigma;0)\subseteq \C^n \to \C$, $y \mapsto f(y)$
analytic and bounded by
$M>0$. 
Then the remainder
\begin{equation} \label{taylor3dar}
h^{(f)}(y) = f(y) - \left( f(0) + \skal{f_y(0)}{y} +\frac{1}{2} \skal{y
  f_{yy}(0)}{y} \right),
\end{equation}
fulfills
for all $\eps \in (0,1)$
the estimate
\begin{equation} \label{taylor3abs}
\left| h^{(f)}(y) \right| \leq \frac{M}{1-\eps} \frac{|y|^3}{\sigma^3}\fa
|y|\leq \eps \sigma.
\end{equation}
\end{lemma}\noindent
{\bf Proof.}
Let us fix an
$\eps$,
$0<\eps < 1$
and 
$y\in \C^n$
with
$|y| \leq \eps \sigma$.
In case
$y=0$
(\ref{taylor3abs})
is an immediate consequence of
(\ref{taylor3dar}).
In case $y$ does not vanish we set
\[ y_0 := \eps \sigma \frac{y}{|y|},\]
such that
$|y_0|=\eps \sigma$, and consider the function
\[ g: \B(\eps^{-1};0 ) \subseteq \C \lto \C, \quad z \mapsto g(z) := f(zy_0).\]
By construction
$g(0) = f(0)$
and with the chain rule we get
\[ g_z (z) = \skal{f_y(z y_0)}{y_0}, \quad g_{zz}(z) = \skal{y_0 f_{yy}(z
  y_0)}{y_0}\fa |z| < \eps^{-1}.\]
Lemma
\ref{taylor3dim1}
yields
\begin{eqnarray*}
\left| h^{(f)}(zy_0) \right| \heq \left| f(zy_0) - f(0) - \skal{f_y(0)}{z y_0}
  -\frac{1}{2} \skal{z y_0
  f_{yy}(0)}{z y_0} \right| \\
\heq \left| g(z) - g(0) - g_z(0) z - \frac{1}{2} g_{zz}(0) z^2 \right| =
\left| h^{(g)} (z) \right|\\
\hleq \frac{M}{1-\eps}\frac{|z|^3}{(\eps^{-1})^3} = \frac{M}{1-\eps} |z|^3
\eps^3\fa |z| \leq \eps (\eps^{-1}) = 1.
\end{eqnarray*}
It is allowed to put 
$z = |y|/(\eps \sigma)$
in this inequality, so
\[ \left| h^{(f)}(zy_0) \right| = \left| h^{(f)} \left( \frac{|y|}{\eps
  \sigma}\,\eps \sigma\frac{y}{|y|} \right) \right| = \left| h^{(f)}(y)
  \right| \leq \frac{M}{1-\eps} \frac{|y|^3}{\eps^3\sigma^3}\eps^3 =
  \frac{M}{1-\eps} \frac{|y|^3}{\sigma^3},\]
and the proof is finished.
\hfill$\Box$
\subsection{Generating symplectic transformations} \label{gct}
\subsubsection*{Auxiliary results on autonomous differential equations}
\begin{theorem} \label{periodischedgl}
Let $\ro>0$, $\Ss(\ro) \subseteq \C^n$, $\V \subseteq \C^m$ open and 
\[f:\Ss(\ro)\times \V \longrightarrow \C^{n+m}, \quad z=(x,y) \mapsto f(z)\]
be continuous and such that
\begin{equation} \label{autonom}
\dot{z} = f(z)
\end{equation}
has unique solutions. The function $f$ shall have the period 
$T>0$
in
$z_1=x_1, \ldots, z_n=x_n$.
We assume that there are numbers
$a, b, \tdelta$, $a\leq0<b$, $0<\tdelta < \ro$
and an open set
$\U\subseteq \V$,
such that the flow
$\fhi$
of 
(\ref{autonom})
exists on
$[a,b) \times \Ss(\ro-\tdelta) \times \U$.
Then the function
\[ \fhi(t,\stv)-\id: \Ss(\ro-\tdelta) \times \U \longrightarrow \Ss(\ro) \times \V, \quad \zeta = (\xi, \eta) \mapsto \fhi(t,\zeta) - \zeta \]
has the period $T$ in
$\zeta_1= \xi_1, \ldots, \zeta_n=\xi_n$
for all $t\in [a,b)$.
\end{theorem}
The assumption on the existence of the flow 
$\fhi$
means, that there is a map
\[ \fhi: [a,b) \times \Ss(\ro-\tdelta) \times \U \lto \Ss(\ro) \times \V\]
with
$\fhi(0,\zeta) = \zeta$
and
$\fhi(\stv,\zeta)$ 
solves the differential equation 
(\ref{autonom}).\\\\
{\bf Proof of theorem \ref{periodischedgl}.}
We show for all
$(t,\zeta) \in [a,b) \times \Ss(\ro-\tdelta) \times \U$
that
\begin{equation} \label{lambada}
\fhi(t,\zeta) +T\cdot e_j = \fhi(t,\zeta+T\cdot e_j) \quad (1\leq j \leq n).
\end{equation}
Let
$j\in \{1,\ldots,n\}$
be arbitrary and set
$h(t) := \fhi(t,\zeta) +T\cdot e_j$ and $g(t) := \fhi(t,\zeta+T\cdot e_j)$.
Then
$h(0) = g(0) = \zeta + T\cdot e_j$
and
\begin{eqnarray*}
\dot{h}(t) \hspace{-2mm}&=&\hspace{-2mm} \dot{\fhi}(t,\zeta) = f(\fhi(t,\zeta)) = f(\fhi(t,\zeta)+T\cdot e_j) = f(h(t)),\\
\dot{g}(t) \hspace{-2mm}&=&\hspace{-2mm} \dot{\fhi} (t, \zeta+T\cdot e_j) = f(\fhi(t,\zeta+T\cdot e_j)) = f(g(t)).
\end{eqnarray*}
Therefore both functions fulfill the differential equation. Hence they coincide. This proves
(\ref{lambada}).
Now
(\ref{lambada})
shows for all
$1\leq j \leq n$
\[ \fhi(t,\zeta+T\cdot e_j) - (\zeta+T\cdot e_j) = \fhi(t, \zeta) - \zeta,\]
which proves the lemma.
\hfill$\Box$
\begin{lemma}
Let
$a<b$ and $f: (a,b) \to \C^m$, $m\in \N$
be an analytic function. Let
$a \leq a_0 < b_0 \leq b$
and suppose that the restriction of $f$ to
$(a_0, b_0)$
maps to 
$\R^m$. Than $f$ maps to
$\R^m$.
\end{lemma}\noindent
{\bf Proof.}
Without loss of generality we may assume
$m=1$,
for in case
$f= (f_1, \ldots, f_m): (a,b) \to \C^m$
is analytic, so is every coordinate function
$f_i$, $1 \leq i \leq m$.
Hence we can apply the lemma for
$m=1$
to each coordinate function and get the result for $f$. So let us assume
$m=1$.\\
Let
$A \subseteq (a,b)$
be the biggest interval, which contains
$(a_0, b_0)$,
and on which $f$ maps to 
$\R^m$. $A$ exists, because it can be constructed as the union of all intervals, which contain
$(a_0, b_0)$
and on which $f$ maps to
$\R^m$.
$A$ is not empty, for it contains
$(a_0, b_0)$.\\\indent
$A$ is closed in
$(a,b)$.
To see that we consider a cluster point
$\alpha$
of 
$A$
and choose a sequence
$\left( x_\ell \right)_{\ell=1}^\infty \subseteq A \setminus \{ \alpha\}$,
which tends to
$\alpha$.
$f$ is in particular continuous on
$(a,b)$, so the limit
\[ f(\alpha) = \lim_{\ell\to \infty} f(x_\ell)\]
exists. It is a limit of real numbers, so it is real as well. Hence
$\alpha \in A$. So $A$ contains its cluster points which means it is closed.\\
However, $A$ is open in
$(a,b)$.
In order to see that consider an arbitrary
$\alpha \in A$.
By assumption $f$ may be expanded in a power series around the point
$\alpha$. The series is given by
\begin{equation} \label{potenzreihe}
f(x) = \sum_{k=0}^\infty \frac{f^{(k)}(\alpha)}{k!} (x-\alpha)^k.
\end{equation}
Herein
$f^{(k)}(\alpha)$
denotes the $k$-th derivative of $f$ in
$\alpha$.
We show that
$f^{(k)}(\alpha)$ is a real number for all $k \in \N_0$.
This is obvious for
$f^{(0)}(\alpha) = f(\alpha)$ because
$\alpha \in A$.
If it is true for some
$k \in \N_0$
then for
$k+1$
as well. Indeed, take a sequence
$\left( x_\ell\right)_{\ell=1}^\infty \subseteq A \setminus \{ \alpha\}$,
which tends to
$\alpha$
and consider the limit
\[ f^{(k+1)}(\alpha) = \lim_{\ell\to\infty} \frac{f^{(k)}(x_\ell) -
  f^{(k)}(\alpha)}{x_\ell-\alpha}.\]
Again, this is a limit of real numbers, hence a real number. So all coefficients of the series
(\ref{potenzreihe}) a real and $f$ maps to
$\R^m$ in a neighborhood of
$\alpha$.
So $\alpha$ is an inner point of $A$ and $A$ is open in
$(a,b)$.\\
Altogether, $A$ is not empty, open and closed in
$(a,b)$, meaning
$A=(a,b)$.
The lemma is proved.
\hfill$\Box$
\begin{theorem} \label{undnochnsatz}
Let
$\ro>0$, $\sigma>0$ and $f\in \Pp_{2n}(\ro,\sigma)$.
Suppose there are
$0 < \tdelta < \ro$, $0 < \eps < \sigma$ and $a \leq 0 < b$
such that the flow
$\fhi$
of the differential equation
\begin{equation} \label{undnochnsystem}
\dot{z} = f(z)
\end{equation}
exists on
$[a,b) \times \D(\ro-\tdelta, \sigma-\eps)$.
If then
$f$ maps real vectors to real values, so does
$\fhi$.
\end{theorem}\noindent
{\bf Proof.}
We consider the restriction of $f$ to real vectors, namely 
\[ g:  \R^n \times \left\{ y\in \R^n\,|\,|y|<\sigma\right\} \lto \R^{2n}, \quad z
\mapsto g(z) := f(z),\]
and the differential equation
\begin{equation} \label{nochnsystem}
\dot{z} =g(z).
\end{equation}
Observe that the domain of $g$ coincides with
$\D(\ro,\sigma) \cap \R^{2n}$.
Now let
\[ \zeta \in \R^n \times \left\{ y\in \R^n\,|\,|y|<\sigma-\eps\right\}\]
be arbitrary. Then there are numbers
$a_1 < 0 <  b_1$ and a solution
\[ h: (a_1, b_1) \lto \R^n \times \left\{ y\in \R^n\,|\,|y|<\sigma\right\} \]
of
(\ref{nochnsystem}). Clearly $h$ is a solution of
(\ref{undnochnsystem})
as well. Therefore
\[ \fhi(t,\zeta) = h(t) \fa t \in (a_1, b_1) \cap [a,b).\]
The set of the $t$ which can applied herein contains an open interval. So the preceding lemma shows that
$\fhi(\stv, \zeta)$
maps to
$\R^{2n}$, which proves the assertion.
\hfill$\Box$
\subsubsection*{Simple canonical transformations}
\begin{theorem} \label{einftrafos1}
Let $\U$, $\V \subseteq \C^n$ be open and connected sets and $Z=(X,Y):\U\times \V \to \C^{2n}$ a simple canonical transformation (see definition \ref{einfachetrafo}). Than we have for all $(\xi,\eta)\in \U\times \V$
\begin{equation} \label{xxi}
\det X_\xi(\xi) \not= 0,
\end{equation}
\begin{equation} \label{yxi}
Y(\xi,\eta) = Y(\xi,0) + \eta X_\xi(\xi)^{-1}.
\end{equation}
\end{theorem}\noindent
{\bf Proof.}
$X$ is independent of $\eta$, so
\[ Z_\zeta = \left( \begin{array}{cc} X_\xi & 0 \\ Y_\xi & Y_\eta \end{array}
\right).\]
Hence
(\ref{defsympl})
implies
\begin{eqnarray*} \maatrix{0}{E_n}{-E_n}{0} \heq \maatrix{X_\xi^{\rm T}}{Y_\xi^{\rm T}}{0}{Y_\eta^{\rm T}}
\maatrix{0}{E_n}{-E_n}{0} \maatrix{X_\xi}{0}{Y_\xi}{Y_\eta} \\
\heq \maatrix{-Y_\xi^{\rm T}}{X_\xi^{\rm T}}{-Y_\eta^\tp}{0}
\maatrix{X_\xi}{0}{Y_\xi}{Y_\eta} = \maatrix{X_\xi^\tp Y_\xi - Y_\xi^\tp
  X_\xi}{X_\xi^\tp Y_\eta}{-Y_\eta^\tp X_\xi}{0}.
\end{eqnarray*}
We consider the right upper block on the left hand and right hand side of the equation and see
\begin{equation} \label{xxiyxi}
X_\xi^\tp Y_\eta = E_n.
\end{equation}
Building determinants we obtain
\[ \det X_\xi (\xi) \det Y_\eta (\xi, \eta) = 1\fa (\xi, \eta) \in \U \times \V.\]
This yields
(\ref{xxi}). Moreover by
(\ref{xxiyxi})
we get
\begin{equation} \label{yetamittp}
Y_\eta = \left( X_\xi^\tp \right)^{-1} = \left( X_\xi^{-1} \right)^\tp.
\end{equation}
Therefore
$Y_\eta$
does not depend on
$\eta$ and consequently
$Y_{\eta\eta} = 0$, such that
$Y$ is affine-linear in
$\eta$. The Taylor expansion of $Y$ with respect to
$\eta$
therefore reads
\[ Y(\xi,\eta) = Y(\xi, 0) + \eta \cdot Y_\eta(\xi, 0)^\tp\fa (\xi, \eta) \in \U \times \V.\]
Together with
(\ref{yetamittp})
we obtain
(\ref{yxi})
and the proof is finished.
\hfill$\Box$
\begin{remark}{\em Theorem
\ref{einftrafos1} 
in particular implies, that simple canonical transformations are affine-linear in $\eta$. So they may be defined for all
$\eta \in \C^n$. 
Moreover the functions
$Y_{k\eta}$ do not depend on $\eta$.
}\end{remark}
Let us denote the uniform convergence of a sequence of functions
$\left( f_k \right)$
on compact subsets of an open set
$\U$
towards some limit function $f$ by
\[ f_k \stackrel{\U,\,{\rm\scriptscriptstyle compact}}{=\!=\!=\!\Longrightarrow} f \qquad (k \to
\infty).\]
Clearly, when
$\U \subseteq \C^n$ or $\U \subseteq \R^n$, the uniform convergence on compact subsets of $\U$ is equivalent to the fact, that the sequence converges uniformly on bounded open subsets of
$\U$.
\begin{theorem} \label{einftrafos2}
Let $\U \subseteq \C^n$ be an open and connected set and
\begin{equation} \label{trafolgar}
Z_k = (X_k, Y_k) : \U \times \C^n \lto \C^n \times \C^n \quad (k \in \N)
\end{equation}
a sequence of simple canonical transformations with the property, that the sequences
$\left( Z_k(\stv,0)\right)_{k=1}^\infty$ and $\left( Y_{k\eta} \right)_{k=1}^\infty$ converge uniformly on compact subsets of $\U$. 
Then
$\left( Z_k \right)_{k=1}^\infty$
converges uniformly on compact subsets of
$\U \times \C^n$ to a simple canonical transformation.
\end{theorem}\noindent
{\bf Proof.}
For all
$k\in \N$
\[ Z_k(\xi, 0) = ( X_k(\xi), Y_k(\xi,0))\]
holds. The functions
$Z_k$
are analytic. By assumption and the theorem of
Weierstrass (see \cite{dieudonne60}, (9.12.1))
there exist analytic functions
$X$, $V$, and $W$,
defined on
$\U$, with
\begin{equation} \label{eifolgen1} 
\komp{\U}{Z_k(\stv, 0)}{(X, V)}\; \mbox{ and }\; \komp{\U}{Y_{k\eta}}{W}, \qquad (k\to \infty).
\end{equation}
The first limit means in particular
\begin{equation} \label{nocheifolgen}
\komp{\U}{X_k}{X}\;\mbox{ and }\; \komp{\U}{Y_k(\stv,0)}{V}, \qquad (k\to \infty).
\end{equation}
By the theorem of Weierstrass we conclude
\[ \komp{\U}{X_{k\xi}}{X_\xi}, \qquad (k \to \infty).\]
Now by
(\ref{yxi})
we have
$Y_{k\eta} = ((X_{k\xi})^{-1})^\tp$,
so the second limit in
(\ref{eifolgen1}) yields for all $\xi \in \U$
\[ E_n = X_{k\xi}(\xi) Y_{k\eta}(\xi)^\tp \lto X_\xi(\xi) W(\xi)^\tp, \qquad
(k \to \infty).\]
Hence
$E_n = X_\xi W^\tp$ holds and
$(X_\xi)^{-1} = W^\tp$
exists, where
\begin{equation} \label{eieifolgen}
\komp{\U}{(X_{k\xi})^{-1}}{(X_\xi)^{-1}}, \qquad (k \to \infty),
\end{equation}
again because of
(\ref{eifolgen1}).
We set
\[ Y(\xi, \eta) := V(\xi) + \eta X_\xi(\xi)^{-1}\fa (\xi, \eta) \in \U
\times \C^n,\]
and show for the functions
$Y_k(\xi, \eta) = Y_k(\xi, 0) + \eta X_{k\xi}(\xi)^{-1}$ that
\begin{equation} \label{polykon}
\kooomp{\U\times \C^n}{Y_k}{Y}, \qquad (k\to \infty).
\end{equation}
For this purpose let
$\mathcal{K}_1 \subseteq \U$ and $\mathcal{K}_2 \subseteq \C^n$
be compact and
$\eps>0$.
By
(\ref{nocheifolgen})
there exists a
$N_1 \in \N$ with
\[ \left| Y_k(\stv, 0 ) - V\right|_{\mathcal{K}_1} < \frac{\eps}{2}\fa k
\geq N_1.\]
Because
$\mathcal{K}_2$
is compact there exists a number
$K>0$, 
such that
$\mathcal{K}_2$
is contained in the ball
$\B(K;0)$.
From
(\ref{eieifolgen})
we infer that there is a
$N_2 \in \N$
with
\[ \left| (X_{k\xi})^{-1} - (X_\xi)^{-1} \right|_{\mathcal{K}_1} <
\frac{\eps}{2nK}\fa k \geq N_2.\]
So for all
$k \geq N_1 + N_2$
\[ \left| Y_k - Y \right|_{\mathcal{K}_1 \times \mathcal{K}_2} \leq \left|
Y_k(\stv, 0 ) - V\right|_{\mathcal{K}_1} + nK \left| (X_{k\xi})^{-1} -
(X_\xi)^{-1} \right|_{\mathcal{K}_1} < \eps\]
holds 
and therefore
(\ref{polykon}) is true. We know from
(\ref{nocheifolgen}) and (\ref{polykon}), that the sequence
$(Z_k)$
converges uniformly on compact subsets of
$\U \times \C^n$
to an analytic function
$Z := (X,Y)$.
It remains to show that $Z$ is a simple canonical transformation. We do already know that $Z$ is analytic and that its component $X$ does not depend on
$\eta$. Hence the only missing information is that $Z$ is a symplectic transformation. Well, by the theorem of Weierstrass we see for all
$(\xi, \eta) \in \U \times \C^n$
\[ J = Z_{k\zeta}(\xi,\eta)^\tp \cdot J \cdot Z_{k\zeta}(\xi,\eta) \lto
Z_{\zeta}(\xi,\eta)^\tp \cdot
J \cdot Z_{\zeta}(\xi,\eta), \qquad (k \to \infty),\]
hence
$Z_\zeta^\tp \cdot J \cdot Z_\zeta = J$.
The proof is finished.
\hfill$\Box$
\subsubsection*{Generating symplectic transformations}
\label{construction.alm}
The discussion in this section is like the one given in
\cite{ruessmann01}. However, we consider an other class of Hamiltonians.
\begin{theorem} \label{satzhamsysteme1}
Let
$K >0$, $\tro>0$, $0<\delta<\tro$ and $0 < \sigma \leq \delta$.
Let 
$F: \D(\tro,\sigma) \to \C$, $F=F(x,y)$
be an analytic function fulfilling
\begin{equation} \label{vorh1}
\left| F_x \right|_{\D(\tro,\sigma)} \leq \frac{K}{\delta}, \quad \left| F_y \right|_{\D(\tro,\sigma)} \leq \frac{K}{\sigma}.
\end{equation}
Then the Hamiltonian system
\begin{equation} \label{kansystem1}
\dot{x} = F_y, \quad \dot{y} = -F_x 
\end{equation}
possesses an analytic flow
\setlength{\mathindent}{1cm}
\[
Z : \left[ 0, \frac{\sigma \delta}{2K} \right) \times \D( \tro-\delta,
  \sigma/2) \lto \D(\tro, \sigma), \quad (t,\zeta) \mapsto Z(t,\zeta),
\]
which is uniquely determined.
\end{theorem}
In particular
$Z(\stv, \zeta)$
is the unique solution to
(\ref{kansystem1})
with respect to the initial value
$Z(0,\zeta) = \zeta \in \D(\tro-\delta, \sigma/2)$. Using the matrix $J$ from definition
\ref{defsymplek}
we can write
(\ref{kansystem1})
in the form
\[ \dot{z} = F_z J^\tp.\]
{\bf Proof of theorem \ref{satzhamsysteme1}.}
The existence theorem of Cauchy
(see \cite{dieudonne60}, (10.4.5))
says, that solutions
$t \mapsto Z(t,\zeta)$
to the initial value
$Z(0, \zeta) = \zeta \in \D(\tro,\sigma)$
exist locally and are uniquely determined. The flow $Z$ is analytic in $t$ and
$\zeta = (\zeta_1, \ldots, \zeta_{2n})$
(see \cite{dieudonne60}, (10.8.2)).
Each solution of
(\ref{kansystem1})
maps to
$\D(\tro,\sigma)$
by definition and it remains to show, that the solutions to the initial values
$\zeta \in \D(\tro-\delta, \sigma/2)$ exist for all $t \in [0, \sigma   \delta /(2K))$.\\
To this end let
$\zeta \in \D(\tro-\delta, \sigma/2)$
be arbitrary. We assume, that the solution
$Z(\stv,\zeta) = (X(\stv,\zeta), Y(\stv,\zeta))$
does only exist up to a
$b\in (0, \sigma\delta /(2 K))$.
By
(\ref{kansystem1})
we have for all
$t\in [0,b)$
\begin{eqnarray*}
X(t,\zeta) -\xi \hspace{-2mm}&=&\hspace{-2mm} \int_0^t F_y(Z(\tau, \zeta))\,d\tau, \\
Y(t,\zeta)-\eta \hspace{-2mm}&=&\hspace{-2mm} \int_0^t -F_x(Z(\tau,\zeta))\,d\tau.
\end{eqnarray*}
Assumption
(\ref{vorh1}) and $0 < b < \sigma \delta / (2 K)$ imply
\begin{eqnarray*}
\left| X(\stv,\zeta)-\xi \right|_{[0,b)} \hspace{-2mm}&\leq&\hspace{-2mm} \sup_{t\in [0,b)} \int_0^t \left|F_y\right|_{\D(\tro,\sigma)}\,d\tau \leq b \frac{K}{\sigma} < \frac{\delta}{2},\\
\left| Y(\stv,\zeta)-\eta\right|_{[0,b)} \hspace{-2mm}&\leq&\hspace{-2mm} \sup_{t\in [0,b)} \int_0^t \left| F_x\right|_{\D(\tro,\sigma)}\,d\tau \leq b\frac{K}{\delta} < \frac{\sigma}{2}.
\end{eqnarray*}
Now let 
$\left(t_k\right)_{k=1}^\infty$ be an increasing sequence in
$[0,b)$
with
$\lim_{k\to \infty} t_k = b$.
According to our assumption on $b$ the sequence
$\left(Z(t_k, \zeta)\right)_{k=1}^\infty$
cannot have a cluster point in
$\D(\tro,\sigma)$
(see \cite{hirsch74}, Chapter 8, \S\,5).
On the other hand, the sequence is contained in the compact set
\[ \left\{ (x,y) \in \C^{2n}\,|\, |\im x|\leq \tro-\delta+b\frac{K}{\sigma},\,|y|\leq \frac{\sigma}{2}+b\frac{K}{\delta} \right\}\subseteq \D(\tro,\sigma),\]
which implies the existence of a cluster point in
$\D(\tro,\sigma)$. This contradiction shows
$b \geq \sigma \delta /(2K)$
and therefore, that the solutions exist for all
$t \in [0, \sigma\delta/(2 K))$.
\hfill$\Box$
\begin{corollary}
\label{korollarhamsysteme1}
Let
$K >0$, $\ro>0$, $0<2\delta<\ro$, and $0 < \sigma \leq \delta$.
Let
$F: \D(\ro,\sigma) \to \C$, $F=F(x,y)$
be analytic and such that
\begin{equation} \label{vorh}
\left| F_x \right|_{\D(\ro,\sigma)} \leq \frac{K}{\delta}, \quad \left| F_y \right|_{\D(\ro,\sigma)} \leq \frac{K}{\sigma}
\end{equation}
holds. Then the Hamiltonian system
\begin{equation} \label{kansystem}
\dot{x} = F_y, \quad \dot{y} = -F_x
\end{equation}
possesses an analytic flow
\setlength{\mathindent}{1cm}
\begin{equation} \label{fluss}
Z : \left[ 0, \frac{\sigma \delta}{2K} \right) \times \D( \ro-2\delta,
  \sigma/2) \lto \D(\ro-\delta, \sigma), \quad (t,\zeta) \mapsto Z(t,\zeta),
\end{equation}
which is uniquely determined.
\end{corollary}\noindent
{\bf Proof.} For the proof it suffices to put
$\tro = \ro-\delta$
in the assumptions of the preceding theorem.
\hfill$\Box$\\\\
When we fix the time $t$ and vary the initial value,
(\ref{fluss}) gives rise to the maps
\begin{equation} \label{hamflussorig}
Z(t,\stv):\D(\ro-2\delta, \sigma/2) \lto \D(\ro-\delta, \sigma), \quad \left( 0 \leq t
< \frac{\sigma \delta}{2K} \right).
\end{equation}
Let us analyze these maps in detail.
\begin{theorem} \label{satzhamsysteme2}
Let
$K >0$, $\ro>0$, $0<2\delta<\ro$, and $0 < \sigma \leq \delta$
with
\[ \frac{\sigma \delta}{2 K} > 1.\]
Let
$F: \D(\ro,\sigma) \to \C$, $F=F(x,y)$
be an analytic function, which is affine-linear in
$y$ and fulfills
(\ref{vorh}).
Then the functions
(\ref{hamflussorig})
satisfy
\begin{equation} \label{flussab1}
\left| Z_\zeta(t,\,\cdot\,)\right|_{\D(\ro-2\delta, \sigma/2)} \leq \exp
\left(\frac{2nK}{\delta \sigma}\,t \right)\fa t \in \left[0, \frac{\sigma\delta}{2 K}\right),
\end{equation} 
\begin{equation} \label{flussab2}
\left| Z_\zeta(t,\stv) -E_{2n}\right|_{\D(\ro-2\delta, \sigma/2)} \leq
\frac{2nK}{\delta \sigma} \exp\left( \frac{2nK}{\delta\sigma}\,t\right)\fa 
t \in [0, 1].
\end{equation}
\end{theorem} \noindent
{\bf Proof.}
We make use of the lemma of Gronwall
(\cite{amann83}, Corollary (6.2)).
For this we have to find an estimate for
$F_{zz}$.
Cauchy's estimate and
(\ref{vorh})
give
\[ \left| F_{xx}\right|_{\D(\ro-\delta,\sigma)} \leq \frac{K}{\delta^2}\leq \frac{K}{\delta \sigma},\quad \left| F_{yx}\right|_{\D(\ro-\delta, \sigma)} \leq \frac{K}{\delta \sigma}.\]
The second inequality and the lemma of Schwarz yield
\[ \left| F_{xy} \right|_{\D(\ro-\delta, \sigma)} = \left| F_{yx}^\tp
\right|_{\D(\ro-\delta, \sigma)} \leq n \left| F_{yx}\right|_{\D(\ro-\delta, \sigma)} \leq \frac{nK}{\delta \sigma}.\]
$F$ is affine-linear in $y$, so
$F_{yy}= 0$. Altogether we obtain
\begin{equation} \label{abshzz}
\left| F_{zz} \right|_{\D(\ro-\delta,\sigma)} \leq (n+1)\frac{K}{\delta \sigma}.\end{equation}
The equation
\[ Z_t^\tp(t,\zeta) = F_z(Z(t,\zeta)) J^\tp\]
holds for all
$0 \leq t < \sigma \delta /(2K)$,
because
$Z(\stv, \zeta)$
solves
(\ref{kansystem}) for all
$\zeta \in \D(\ro-2\delta, \sigma/2)$.
(On the left hand side we have to write
$Z_t^\tp$
because of our definition
$\dot{Z} = Z_t^\tp$ on page
\pageref{zeitabl}.)
Differentiating with respect to 
$\zeta$
yields
\begin{equation} \label{zetawitz}
Z_{\zeta t}(t,\zeta) = (Z_t^\tp)_\zeta(t, \zeta) = J F_{zz}(Z(t,\zeta))
\cdot Z_\zeta(t,\zeta).
\end{equation}
Now integration with respect to $t$ gives
\begin{equation} \label{dirkiint}
Z_\zeta(t,\zeta) = E_{2n} + \int_0^t J F_{zz}(Z(\tau,\zeta))\cdot
Z_\zeta(\tau,\zeta)\,d\tau.
\end{equation}
With
(\ref{abshzz})
we obtain the estimate
\begin{eqnarray*}
\lefteqn{|Z_\zeta (t,\zeta)| \leq 1 + \int_0^t \left| F_{zz}\right|_{\D(\ro-\delta, \sigma)} |Z_\zeta(\tau, \zeta)|\,d\tau \leq 1 + \frac{(n+1)K}{\delta \sigma} \int_0^t |Z_\zeta(\tau, \zeta)|\,d\tau}\\
&&\hspace{6.2cm} \fa \zeta \in \D(\ro-2\delta, \sigma/2),\, t \in \left[0, \frac{\sigma\delta}{2K}\right).
\end{eqnarray*}
With the lemma of Gronwall
\begin{eqnarray*}
\lefteqn{|Z_\zeta(t,\zeta)| \leq \exp \left(\frac{(n+1)K}{\delta \sigma}\,t
  \right) < \exp \left( \frac{2nK}{\delta \sigma}\, t
  \right)}\\&&\hspace{6.2cm} \fa \zeta \in \D(\ro-2\delta, \sigma/2),\, t \in
\left[0, \frac{\sigma\delta}{2K}\right)
\end{eqnarray*}
follows. To obtain the second estimate, we derive with
(\ref{dirkiint}) for $t \in [0, \sigma \delta /(2K))$ and $\zeta \in \D(\ro-2\delta, \sigma/2)$
\[ Z_\zeta(t,\zeta) - E_{2n} = \int_0^t JF_{zz}(Z(\tau, \zeta))\,d\tau + \int_0^t JF_{zz}(Z(\tau, \zeta))(Z_\zeta(\tau, \zeta)-E_{2n})\,d\tau.\]
This together with
(\ref{abshzz})
implies 
\begin{eqnarray*}
\left| Z_\zeta(t,\zeta)-E_{2n}\right| \hleq \frac{(n+1)K}{\delta \sigma}\, t + \frac{(n+1)K}{\delta\sigma} \int_0^t \left| Z_\zeta(\tau,\zeta)-E_{2n}\right|\,d\tau\\ \hleq \frac{2nK}{\delta\sigma} + \frac{2nK}{\delta\sigma} \int_0^t \left| Z_\zeta(\tau,\zeta)-E_{2n}\right|\,d\tau\\&&\hspace{4.2cm} \fa \zeta \in \D(\ro-2\delta, \sigma/2), \,t\in [0,1].
\end{eqnarray*}
Here the lemma of Gronwall says
\[\left| Z_\zeta(t,\zeta)-E_{2n}\right| \leq \frac{2nK}{\delta\sigma} \exp \left(\frac{2nK}{\delta\sigma} \,t\right)\fa \zeta \in \D(\ro-2\delta, \sigma/2),\, t\in [0,1].\]
The theorem is shown.
\hfill$\Box$
\begin{theorem} \label{satz2trafos1}
Let
$K >0$, $\ro>0$, $0<2\delta<\ro$, and $0 < \sigma \leq \delta$.
Let the function
$F: \D(\ro,\sigma) \to \C$, $F=F(x,y)$
be analytic and such that
(\ref{vorh}) holds. Then the maps 
 (\ref{hamflussorig})
are symplectic transformations.
\end{theorem} \noindent
{\bf Proof.}
We meet the assumptions of corollary
\ref{korollarhamsysteme1}. Therefore the flow
(\ref{fluss})
and the maps
(\ref{hamflussorig})
exist. We have to prove:
\begin{equation} \label{telefon}
Z_\zeta(t,\zeta)^\tp J Z_\zeta(t,\zeta)= J\fa (t,\zeta) \in \left[ 0,
\frac{\sigma\delta}{2K}\right) \times \D(\ro-2\delta, \sigma/2).
\end{equation}
This equation is certainly true for
$t=0$, because
$Z(0,\stv)$
is the identity and so
$Z_\zeta(0, \zeta) = E_{2n}$ for all $\zeta \in \D(\ro-2\delta, \sigma/2)$.\\
To get the assertion tor all
$t \in [0, \sigma \delta /(2K))$
we show that the left hand side of
(\ref{telefon})
is constant with respect to $t$. To this end we calculate for
$(t,\zeta) \in [0, \sigma \delta/(2K))\times \D(\ro-2\delta, \sigma/2)$
with
(\ref{zetawitz})
\begin{eqnarray*}
\lefteqn{\frac{\partial}{\partial t} \left( Z_\zeta(t,\zeta)^\tp J Z_\zeta(t,\zeta) \right) =
 (Z_\zeta^\tp)_t(t,\zeta) J Z_\zeta(t,\zeta) + Z_\zeta(t,\zeta)^\tp J Z_{\zeta
  t}(t,\zeta)} \\ &&= Z_\zeta(t,\zeta)^\tp F_{zz}(Z(t,\zeta))J^\tp J
 Z_\zeta(t,\zeta) + Z_\zeta(t,\zeta)^\tp JJF_{zz}(Z(t,\zeta))Z_\zeta(t,\zeta) \\&\hspace{1cm}&= Z_\zeta(t,\zeta)^\tp F_{zz}(Z(t,\zeta))
 Z_\zeta(t,\zeta) - Z_\zeta(t,\zeta)^\tp F_{zz}(Z(t,\zeta))Z_\zeta(t,\zeta) =
 0.\end{eqnarray*}
This ends the proof.
\hfill$\Box$
\begin{theorem}\label{haffinlinear}
Let
$K >0$, $\ro>0$, $0<2\delta<\ro$, and $0 < \sigma \leq \delta$.
Let
$F: \D(\ro,\sigma) \to \C$, $F=F(x,y)$
be analytic, so that
(\ref{vorh}) holds, and affine-linear in $y$. Then the maps
(\ref{hamflussorig})
are simple canonical transformations.
\end{theorem}\noindent
{\bf Proof.} 
The assumptions on the function $F$ mean, that $F$ can be written as
\[ F(x,y) = F_1(x) + \skal{y}{F_2(x)},\]
where
$F_1:\Ss(\ro)\to \C$ and $F_2:\Ss(\ro) \to \C^n$
are analytic functions. System
(\ref{kansystem})
reads in this case
\[ \dot{x} = F_2(x), \quad \dot{y} = -F_{1x}(x) - y\cdot F_{2x}(x).\]
The first equation possesses a unique solution
$\widetilde{X}(\stv,\xi)$, $\widetilde{X}(0,\xi) = \xi$
for all initial values
$\xi \in \Ss(\ro-2\delta)$.
This solution exists for all
$0 \leq t < \sigma\delta /(2 K)$,
as can be seen as above. Let us consider the system
\begin{equation} \label{kannsystem}
\dot{x} = F_2(x), \quad \dot{y} = 0.
\end{equation}
Obviously its solutions are given by
$\widetilde{Z} (\stv,\xi,\eta) = (\widetilde{X}(\stv,\xi),\eta)$.
Now, let
$Z=(X,Y)$
be a solution of
(\ref{kansystem})
with initial value
$Z(0,\zeta) = \zeta = (\xi,\eta)$.
Then
$X(0,\zeta) = \xi$
holds and
$t \mapsto (X(t,\zeta), \eta)$
solves
(\ref{kannsystem}).
Therefore $X$ has the same values as
$\widetilde{X}$, meaning
\setlength{\mathindent}{1cm}
\[ X(t,\xi,\eta) = \widetilde{X}(t,\xi)\fa (t,\xi,\eta) \in \left[ 0,
  \frac{\sigma \delta}{2K} \right) \times \D(\ro-2\delta, \sigma/2).\]
Hence $X$ is independent of
$\eta$
and the map
(\ref{hamflussorig})
is a simple canonical transformation as was to be shown.
\hfill$\Box$\\\\
We resume the results of this appendix
\ref{gct}
in the following theorem.
\begin{theorem} \label{trafoterminator}
Let
$K >0$, $\ro>0$, $0<2\delta<\ro$, and $0 < \sigma \leq \delta$ with
\[ \frac{\sigma\delta}{2K}> 1.\]
Let the function $F \in \Pp(\ro,\sigma)$ 
fulfill estimates
(\ref{vorh}) and be affine-linear in $y$. 
Then the maps
(\ref{hamflussorig})
are simple canonical transformations, for all $0 \leq
 t < \sigma \delta /(2K)$ we have $Z(t,\stv) -\id \in \Pp_{2n}(\ro-2\delta, \sigma/2)$,
 and the estimates (\ref{flussab1}) and (\ref{flussab2}) are fulfilled.
\end{theorem}\noindent
{\bf Proof.}
The maps
(\ref{hamflussorig}) are well-defined and analytic by corollary
\ref{korollarhamsysteme1}. 
They are simple canonical transformations by theorem
\ref{haffinlinear}. 
The assumptions of
theorem \ref{periodischedgl} are met, one has to put
\[ \V = \B(\sigma;0),\, f = F_z J^\tp, \,T = 2\pi, \,a=0, \,b= \sigma
\delta/(2K),\]\[\tdelta = 2\delta,
\,\U = \B(\sigma/2;0) \mbox{ and } \fhi = Z.\]
Therefore
$Z(t,\stv) - \id$
has period $2\pi$ in $x$ for all $0 \leq t < \sigma\delta/(2K)$.
The assumptions of theorem \ref{undnochnsatz}
are achieved with
\[ f=F_z J^\tp,\, \tdelta = 2\delta,\, \eps = \sigma/2, \,a=0, \,b = \sigma \delta/(2K) \mbox{ and }\fhi =
Z.\]
So $Z$ maps real vectors to real values.
This shows $Z(t,\stv) -\id \in \Pp_{2n}(\ro-2\delta, \sigma/2)$.
Finally
(\ref{flussab1}) and (\ref{flussab2}) are a consequence of theorem
\ref{satzhamsysteme2}. This finishes the proof.\\${}$\hfill$\Box$
\end{appendix}

\end{document}